\documentclass[10pt]{amsart}
\usepackage{amsfonts, amssymb, amsmath, latexsym}

\newcommand{\la}{\langle}
\newcommand{\ra}{\rangle}

\newcommand{\goth}{\mathfrak}

\newcommand{\C}{\mathbb{C}}

\newcommand{\R}{\mathbb{R}}

\newcommand{\bz}{\bar{\partial}}

\newcommand{\Pz}{\partial_z}
\newcommand{\Bz}{\partial_{\bar{z}}}

\newcommand{\bb}{\mathbb} 
\newcommand{\II}{I\! \!I}

\newcommand{\h}{\overset{\rightarrow}{H}}
\newtheorem{theorem}{Theorem}[section]
\newtheorem{proposition}{Proposition}[section]

\newtheorem{lemma}[theorem]{Lemma}
\newtheorem{remark}{Remark}[section]

\begin{document}

\title{ conformal geometry of marginally trapped surfaces  in  $\bb{S}^4_1$ } 
\author{E. Hulett}

 \thanks{Partially supported by research grants from CONICET, SECYT-UNC  and 
FONCyT Argentina.} 
\address{C.I.E.M. - Fa.M.A.F.
 Universidad Nacional de C\'ordoba,
Ciudad Universitaria,
5000 C\'ordoba, Argentina.
phone/fax: +54 351 4334052/51}
\email{ hulett@famaf.unc.edu.ar}
\date{}

 \begin{abstract}
A spacelike surface $S\subset \bb{S}^4_1$ is marginally trapped if its mean curvature vector is lightlike.  On any oriented spacelike surface $S \subset \bb{S}^4_1$ we show that a choice of orientation of the  normal bundle $\nu(S)$        determines a smooth map $G: S \to \bb{S}^3$ which we call the null Gauss map of $S$.  We show that if $S$ is marginally trapped then  $G$ is a conformal immersion away  the zeros of certain  quadratic Hopf-differential of $S$ and so the surface  $G(S)$ is uniquely determined up to conformal  transformations of $\bb{S}^3$ by  two  invariants: the normal Hopf differential $\kappa$ and the Schwartzian derivative $s$. We show that these invariants plus an additional  quadratic differential $\delta$ are related by a differential equation and determine the geometry of $S$ up to ambient isometries of $\bb{S}^4_1$.  This allows us to  obtain a   characterization of     marginally trapped surfaces $S$ whose null Gauss image is a  {\it constrained Willmore}  surface in $\bb{S}^3$~\cite{bohle}.  As an application of these results  we construct and study    integrable  non-trivial one-parameter  deformations of  marginally trapped surfaces with non-zero parallel mean curvature vector and those   with flat normal bundle. 
 \end{abstract}

\maketitle

\section{introduction}
   
A spacelike surface  immersed in a $4$-dimensional Lorentz manifold   is called marginally trapped if its  mean curvature vector is  everywhere null or lightlike. 
The notion   of marginally trapped surfaces  was introduced by R. Penrose  and  plays a key role  in   the  sigularity theory    of  Einstein's equations~\cite{bychen}.  
The marginally trapped equation $\la \h, \h \ra =0$ is interpreted in relativity theory as the condition describing the event horizon of a black hole~\cite{bychen}. 
In    differential geometry   marginally trapped surfaces are  viewed as    
natural  generalizations of minimal   surfaces.  \\
Different  aspects of the geometry of marginally trapped surfaces have drawn the attention of geometers recently. In~\cite{aledo-galvez-mira},~\cite{Liu} the authors provide different Weierstrass-type representation formulas of marginally trapped surfaces in $\R^4_1$. Also~\cite{chen-vanderveken} and~\cite{bychen} deal with a  classifications of marginally trapped surfaces with parallel mean curvature. 
The notion of marginally trappedness has also been  considered recently in  higher dimensions and co-dimensions  with  very interesting results, see~\cite{anciaux},~\cite{anciaux-godoy}. \\

		Our   goal in  this paper   is to study  geometric properties of  oriented marginally trapped surfaces in $\bb{S}^4_1$  in terms of   conformal invariants of these surfaces.  As  applications  of these ideas we construct  spectral  non-trivial deformations of marginally trapped  surfaces with parallel non-zero mean curvature and those with flat normal bundle.  \\
More specifically,  given  an  oriented spacelike surface  $S \subset \bb{S}^4_1$ its normal bundle $\nu(S)$ is Lorentzian hence   at each point $x$ of $S$ there are  two linearly independent  null directions  say, $n_+(x), n_-(x)$ which vary smoothly  with $x$ and  determine a pair of  smooth maps from $S$  to the $3$-sphere $\bb{S}^3$, viewed as the manifold of null directions of Minkowski space $\R^5_1$.  Such    maps can be  interpreted  as    (pseudo) inverses of the conformal Gauss map $Y$ introduced by R. L. Bryant~\cite{bryant}. We show that a  choice of  orientation on $\nu(S)$ distinguishes a preferred map, say $n_+$   which we call the {\it null Gauss map} $G$ of the spacelike surface $S$. 
	When  $S$ is marginally trapped and certain Hopf quadratic differential is never zero on $S$,   then $G$ is a conformal immersion of the surface $S$ into the conformal sphere $\bb{S}^3$ and  its geometry is dictated by two  conformal  invariants: the Schwartzian $s$ and the normal Hopf differential $\kappa$. These invariants  were  introduced and studied   by Burstall et al. in~\cite{burstall-pedit-pinkall}, see also~\cite{xiangma-thesis}.  \\

In Section~\ref{KonformalinvariantsMTS} we   obtain an equation which relates   the
  conformal invariants $s ,\kappa$   with the $\delta$-quadratic differential,  a new geometric invariant    of the corresponding  marginally trapped surface.  
As a first     consequence of  this equation we prove     Theorem~\ref{GaussConstWillmore}  which says that the null Gauss map of an oriented  marginally trapped surface $S$ is a  constrained Willmore  surface in $\bb{S}^3$ if and only if $S$  has non-zero parallel mean curvature vector.
Constrained Willmore surfaces  
were introduced and studied in~\cite{bohle}. They   are defined as extremes of the Willmore energy with respect to variations preserving the underlying conformal structure of the surface. 
A second  consequence is Theorem~\ref{KSdetermineDelta} which states that a marginally trapped surface is essentially determined up to ambient isometries by the conformal invariants $\kappa,s$ of its null Gauss map.\\
In Section~\ref{preliminaries}  we fix notations and derive the structure equations of spacelike surfaces in $\bb{S}^4_1$.    Section~\ref{conformalsurfacetheoryS3} contains  a short survey  of $O(3,1)$-invariant geometry of surfaces in the conformal sphere $\bb{S}^3$. For a detailed exposition on  the conformal invariant  geometry of surfaces in $\bb{S}^n$ see~\cite{burstall-pedit-pinkall},~\cite{jeromin} and~\cite{xiangma-thesis}. \\
In Section~\ref{deformationsMTS} we consider  marginally trapped surfaces admiting non-trivial integrable one parameter deformations.  The  deformation is induced by  a spectral parameter which determine  symmetries of the compatibility equations thus   giving rise to one-parameter families of surfaces obtained by deformation of a given surface. We consider  here  deformations of    two   kinds  of   marginally trapped surfaces in $\bb{S}^4_1$ namely,  
  surfaces  with non-zero parallel mean curvature vector, and surfaces   with flat normal bundle.  
			 In the first  case we show that the deformation originates in  the associated family of an auxiliar  harmonic map $\phi$ with values in a pseudo riemannian complex quadric $Q$.  We show that a marginally trapped surface $f: \Sigma \to \bb{S}^4_1$ has constrained Willmore null Gauss map if an only if  an  auxiliar   map $\phi$ with values in the complex quadric $Q$ is harmonic. We use the associated family of $\phi$ to obtain  a  symmetry  of the  compatibility equations of the constrained Willmore null Gauss map of $f$, thus giving rise to the associated family of  $f$ and its null Gauss map.  
		For marginally trapped surfaces with flat normal bundle  we obtain a one-parameter  deformation which originates in   the so-called {\it Calapso-Bianchi} or isothermic T-tranformation of isothermic surfaces in $\bb{S}^3$~\cite{burstall-pedit-pinkall},~\cite{xiangma-thesis}.
	Motivated by~\cite{burstall-pedit-pinkall} we show that  both deformations may  be unified in  an extended action of $\C-\{ 0 \}$ on  the class of marginally trapped surfaces with non-zero parallel mean curvature.  
				We conclude with a description of this  extended action on non-isotropic marginally trapped tori with non-zero parallel mean curvature vector.\\

\section{preliminaries} \label{preliminaries}

Denote by $\R^5_1$ the real $5$-dimensional vector space    with canonical coordinates $(x_0, x_1, x_2, x_3, x_4)$   equipped with the Lorentz inner product 
\begin{equation}\label{lorentzmetricR51}
 \la x,y \ra = x_0 y_0 + x_1 y_1 + x_2 y_2 + x_3 y_3 - x_4 y_4.
\end{equation} 
 De Sitter $4$-space    is defined as the unit sphere in $\R^5_1$: 
$$\bb{S}^4_1 = \{ x \in \R^5_1 : \la x,x \ra =1 \}$$
Thus  $\bb{S}^4_1$ is a connected simply connected $4$-dimensional manifold which  inherits from  $\R^5_1$  a lorentzian metric $\la .,. \ra$ of constant sectional curvature $+1$.  
The    complex bilinear  extension of   the Lorentz metric to $\C^5$ is given by 
$$  \la z,w \ra =  z_0 w_0 + z_1 w_1 + z_2 w_2 + z_3 w_3 - z_4 w_4$$
 and the corresponding (pseudo) hermitian inner product  is given   by $\la z,\bar{w} \ra$. We denote by $\C^5_1$  the complex space $\C^5$ endowed with the inner product $\la z,\bar{w} \ra$.\\
The Lie group $SO(4,1)$ acts transitively on $\bb{S}^4_1$ by isometries, so that choosing $e_0 \in \bb{S}^4_1$ as the base point, then $\bb{S}^4_1$ is isometric to the (pseudo) riemmanian symmetric space $SO(4,1) /SO(3,1)$. \\  		
 A  non-zero   vector   $X \in \R^5_1$ is said to be {\it future pointing}  if  $\la X, e_4 \ra < 0$.  This induces  a time orientation on $\bb{S}^4_1$:   a non-zero tangent vector  $X \in T_p \bb{S}^4_1$ is  future pointing if its tranlated    to the origin is future pointing.  
If  $X$ is future pointing  and satisfies $\la X, X \ra =-1$, then  (its tranlated) $X$ lies in the  real $4$-hyperbolic space 
			$  \bb{H}^4 = \{ x \in \R^5_1: \la x,x\ra =-1, x_4 >0 \}$. \\ 
	
 Let $\Sigma$ be a connected orientable surface and $f : \Sigma \to \bb{S}^4_1$ a spacelike  immersion i.e.  the induced metric $g= f^{*} \la .,.\ra$ is Riemannian and it   determines a conformal structure on $\Sigma$. Then    $f$   preserves this  conformal structure i.e. 
   $\la 
  f_z ,  f_z \ra=0$,  for every  local complex
  coordinate $z=x+iy$ on $\Sigma$, where 
  $\Pz = \frac{1}{2} ( \frac{\partial}{\partial x}- i \frac{\partial}{\partial y})$, and $\Bz = \frac{1}{2} ( \frac{\partial}{\partial x} + i \frac{\partial}{\partial y})$, 
are the  complex partial operators.
    Equivalently,  
   \begin{equation} \label{confID}
       \la f_x , f_y \ra =0, \quad \| f_x \|^2 =\| f_y \|^2 >0.
  \end{equation}
Conversely, if $f: \Sigma \to \bb{S}^4_1$ is a conformal immersion from a Riemann surface, then $\la f_x , f_y \ra =0$, and $ \| f_x \|^2 =\| f_y \|^2 \neq 0$, for every local complex coordinate $z=x+iy$. Since the ambient $\bb{S}^4_1$ is lorentzian,  $f_x, f_y$ have positive  squared norm $\| f_x \|^2 =\| f_y \|^2 >0$, and so   $ f: (\Sigma,g) \to \bb{S}^4_1$ is a spacelike isometric immersion, where $g$ is the induced metric. Respect to a  local complex coordinate  $z=x+iy$ on $\Sigma$ we introduce 
          a conformal parameter $u$ by $\la  f_z, f_{\bar{z}} \ra=e^{2u}$, so that    $g=2 e^{2u}( dx^2 +   dy^2)$ is the local expression of the induced metric. 	
        Since $f$ is conformal we have  $2\la  f_{\bar{z} z} ,  f_z \ra = \Pz \la  f_z,  f_z
    \ra =0$ and also $2\la  f_{\bar{z} z} , f_{\bar{z} } \ra = \Pz \la f_{\bar{z} }, f_{\bar{z} }
    \ra =0$, which says that     $f_{\bar{z} z}$ has no tangential component.\\
The second fundamental form  of $f$ is defined by  $\la \II(X,Y), N \ra = - \la df(X), dN(Y) \ra$,  for $ X, Y \in T \Sigma$ and for every normal field $N$ along $f$.  The mean curvature vector of $f$ is the trace of $\II$:  $\h:= \frac{1}{2}trace \II $. Since $f_{z \bar{z}}$ has no tangential component it decomposes into its $f$ and $\h$ components by      $ f_{z \bar{z} } = -e^{2u}f + e^{2u} \h$.   \\

The pullback bundle of the tangent bundle of $\bb{S}^4_1$  decomposes into the tangent bundle and the normal bundle of $f$: 
$f^{*}(T \bb{S}^4_1) = T \Sigma \oplus \nu(f)$.  
Since  $f: \Sigma \to \bb{S}^4_1$ is spacelike and $\Sigma$ is orientable,   the  normal bundle $\nu(f)$  is an  orientable lorentzian vector bundle. Fixing   an orientation on $\nu(f)$, let $\{ N_1, N_2 \} \subset \Gamma(\nu(f))$ be  an (ordered) orthonormal frame  satisfying  
$$
\la N_2, N_2 \ra =-1, \quad  \la N_1, N_2 \ra =0,  \quad \la N_1, N_1\ra =1.
$$
If  we demand that       $N_2$ be  future pointing,  then either $\{ N_1, N_2 \}$ has the same orientation as $\nu(f)$, or $\{ -N_1, N_2 \}$ has the same orientation as $\nu(f)$.   We  say that  an    orthonormal frame  $\{ N_1, N_2 \} \subset \Gamma(\nu(f))$ is {\it positively oriented}, $\{ N_1, N_2 \}$ has the same orientation as $\nu(f)$ and  $N_2$ is    (timelike)  future pointing. Note that if $\{ N_1, N_2 \}$ is positively oriented then $\{ -N_1, -N_2 \}$ has the same orientation as $\nu(f)$, but it is not positively oriented since $-N_2$ points to the past. \\   		 
			
In terms of a normal  orthonormal frame the second fundamental form is given by
    \begin{equation}\label{2ndfundform}
    \II =   -\la df, dN_1\ra N_1 + \la df, dN_2 \ra N_2.
    \end{equation}
Let    $\xi_1 := \la f_{zz}, N_1\ra$, $\xi_2 := -\la f_{zz}, N_2 \ra$. Since $f$ is conformal an easy calculation gives  $f_{zz} = 2 u_z f_z+ \xi_1 N_1 + \xi_2 N_2$. In particular  the $(2,0)$-part of $\II$ is given by
  \begin{equation}\label{hol2ndfundform}
  \II( \Pz, \Pz) = \xi_1 N_1 + \xi_2 N_2
  \end{equation}   
 Set   $h_1 := \la \h, N_1 \ra$ and  $h_2 := - \la \h, N_2 \ra$, then    $\h = h_1 N_1 + h_2 N_2$.\\
The structure equations of a conformal immersion $f: \Sigma \to \bb{S}^4_1$ are given by  
\begin{equation}\label{structureEqs2}
    \begin{array}{l}
    f_{zz} =  2 u_z f_z+ \xi_1 N_1 + \xi_2 N_2 \\
    f_{\bar{z}z} =  -e^{2u}f + e^{2u} \h,\\
    \Pz N_1 = -h_1 f_z - e^{-2u} \xi_1 f_{\bar{z}} + \sigma N_2,\\
    \Pz N_2 = h_2 f_z + e^{-2u} \xi_2 f_{\bar{z}} + \sigma N_1,
    \end{array}
    \end{equation}
and the compatibility among these equations is just Gauss's, Codazzi's and Ricci's equation:
    \begin{equation}\label{gausscodazziricciBasic}
    \begin{array}{ll}
    \text{\sf Gauss,} &  2 u_{\bar{z}z} = -e^{2u} +e^{-2u} (|\xi_1|^2 - |\xi_2|^2) - e^{2u} \| \h \|^2,\\
    \text{\sf Codazzi,} &    e^{2u}(\Pz h_1 + \sigma h_2 )= \Bz \xi_1 + \xi_2 \bar{\sigma},\\
		   & e^{2u}(\Pz h_2 + \sigma h_1 ) = \Bz \xi_2 +  \xi_1 \bar{\sigma},\\ 
    \text{\sf Ricci,} &  Im (  \sigma_{\bar{z}} ) = e^{-2u} Im (\xi_1 \bar{\xi}_2).
    \end{array}
    \end{equation}

A spacelike surface $f: \Sigma \to \bb{S}^4_1$ is called { \it marginally trapped}  if its mean curvature vector is  null or lightlike:  $\la \h, \h \ra =0$.  
	If     $\h \neq 0$, then after a change of orientation of 
the normal bundle (i.e. after a change of sign $N_1 \mapsto -N_1$) if necessary,   
	the   marginally trapped condition $ \la \h ,\h \ra =h_1^2-h_2^2 = 0$,   reads  $h_1 = h_2$,   with  $h_1 + h_2 \neq 0$, in this case  the mean curvature vector satisfies      
\begin{equation}\label{MCV}
\h = h(N_1+N_2), \quad \text{with} \,\, h=h_1=h_2.
\end{equation}
We call  $h$ the  {\it the mean curvature function} of $f$ respect to the  positively oriented lorentzian normal frame $\{ N_1, N_2 \}$. 
From the second structure equation 	$f_{\bar{z}z} =  -e^{2u}f + e^{2u} \h$, it follows that  the mean curvature function $h$ satisfies  
\begin{equation}\label{meancurvature}
h= e^{-2u} \la f_{\bar{z}z}, N_1  \ra = - e^{-2u} \la f_{\bar{z}z}, N_2  \ra. 
\end{equation}
Hence   $f$ is marginally trapped iif 
$\la f_{\bar{z}z}, (N_1 + N_2) \ra =0$. 
	Introducing the  function $\sigma :=  - \la \Pz N_1 , N_2 \ra$,  the structure equations of a marginally trapped immersion $f$ are  thus given by 
        \begin{equation}\label{structureEqs}
    \begin{array}{l}
    f_{zz} =  2 u_z f_z+ \xi_1 N_1 + \xi_2 N_2, \\
    f_{\bar{z}z} =  -e^{2u}f + e^{2u} \h,\\
    \Pz N_1 = -h f_z - e^{-2u} \xi_1 f_{\bar{z}} + \sigma N_2,\\
    \Pz N_2 = h f_z + e^{-2u} \xi_2 f_{\bar{z}} + \sigma N_1,
    \end{array}
    \end{equation}
where $\{ N_1, N_2 \} \subset \Gamma(\nu(f))$ is  a positively oriented orthonormal frame. 
    The compatibility conditions or equations of Gauss, Codazzi and Ricci  above reduce to:  
    \begin{equation}\label{gausscodazziricci}
    \begin{array}{ll}
      2 u_{\bar{z}z} = -e^{2u} +e^{-2u} (|\xi_1|^2 - |\xi_2|^2),\\
         e^{-2u}( \xi_{1 \bar{z}} + \xi_2 \bar{\sigma}) = ( h_z + \sigma h ), \\
		 e^{-2u}( \xi_{2 \bar{z}} +  \xi_1 \bar{\sigma})= (h_z + \sigma h ),\\ 
      Im (  \sigma_{\bar{z}} ) = e^{-2u} Im (\xi_1 \bar{\xi}_2).
    \end{array}
    \end{equation}

 The  Gaussian curvature of the induced metric $g$   is given by 
$K= - \Delta_g u = -2 e^{-2u}  u_{\bar{z} z}$,
 where $\Delta_g=  2 e^{-2u} \Bz \Pz$,  is the Laplace operator of the induced  metric $g$.  From  Gauss equation~\eqref{gausscodazziricci} we obtain  the expression of the Gaussian curvature of the induced metric on $\Sigma$,
 \begin{equation}\label{gausscurvature} 
K =  1   - e^{-4u} (|\xi_1|^2 - |\xi_2|^2),
\end{equation}

Let    $\nabla^{\bot}$ denote the covariant derivative on  the normal bundle $\nu(f)$, then  $\omega := \la \nabla^{\bot} N_2 , N_1 \ra$ is the corresponding connection one form.   Fixed an orientation on the normal bundle $\nu(f)$  the normal curvature is defined by $d \omega = K^{\bot} dA_g$, where $dA_g$ is the area form of the induced metric $g$.  Thus   $\omega = 2 Re (\sigma dz)$,  and so  
$ d \omega  =  -4 Im (\sigma_{\bar{z}}) dx \wedge dy$.    From Ricci equation above it follows that   $ Im (  \sigma_{\bar{z}} ) = e^{-2u} Im (\xi_1 \bar{\xi}_2) $.  Thus  since  $dA_g = 2 e^{2u} dx \wedge dy$,  the normal curvature function is given by 
\begin{equation}\label{normalcurvature}
K^{\bot}= -e^{-2u} Im (\sigma_{\bar{z}})= -e^{-2u}Im (\xi_1 \bar{\xi}_2).
\end{equation}
Then the normal bundle is flat if and only if $K^{\bot} =0$. \\

On the other hand from Codazzi's equation  the covariant derivative of the mean curvature vector of a conformal immersion $f: \Sigma \to \bb{S}^4_1$ is given by  
\begin{equation}\label{nablaH}
\nabla^{\bot}_{\Pz} \h = (\Pz h_1 + \sigma h_2 ) N_1 + (\Pz h_2 + \sigma h_1 ) N_2.    
\end{equation}
 In particular if  $f$ is marginally trapped  then in  a positively oriented orthonormal frame $\{ N_1, N_2 \} \subset \Gamma(\nu(f))$  the above formula becomes 
\begin{equation}\label{normalderivativeH}
 \nabla^{\bot}_{\Pz} \h = e^{-2u} ( h_z + \sigma h ) (N_1 +  N_2).
\end{equation}
Hence    $f$ has parallel mean curvature vector if and only if  $h_z + \sigma h =0$.

\section{  surface theory  in the conformal sphere $\bb{S}^3$}\label{conformalsurfacetheoryS3}
We give a brief account  of    Moebius surface geometry   in $\bb{S}^3$  such as  exposed in~\cite{burstall-pedit-pinkall}. For  detailed proofs  and further developements  we refer the reader to~\cite{burstall-pedit-pinkall},~\cite{jeromin} and~\cite{xiangma-thesis}.\\

  The null or light cone in $\R^5_1$ is defined by       
\begin{equation}\label{nullconeR5} 
 \mathcal{L}= \{ 0 \neq x \in \R^5_1 : \la x,x \ra =0 \}.
\end{equation}
The {\it future light cone} $\mathcal{L}_+ \subset \mathcal{L}$ consists of  future pointing   vectors  $ x \in \mathcal{L}$.
For every  $x \in \bb{S}^3 \subset \R^4$,  the point  $(x, 1)\in \R^5_1$ lies  in the future light cone $\mathcal{L}_+$.  We are using here the fact that any  vector $x$ in $\R^5_1$ may be uniquely written as  an ordered pair $(x', t)$ with $x' \in \R^4$ and $t \in \R$, thus giving rise to an isomorphism  $\R^4 \oplus \R \to   \R^5_1$. In particular  points on $\mathcal{L}$  are of the form $(x, \pm \|x \|^2)$, with $x \in \R^4$.\\  
The  map  $ \bb{S}^3 \ni x  \mapsto [(x, 1)]$ identifies the unit round sphere $\bb{S}^3 \subset \R^4$ with the projectivization of the light cone, $P(\mathcal{L}) \subset \bb{RP}^5$. Let   $O_+(4,1)$ be the group  of orthogonal transformations of $\R^5_1$ preserving the time orientation. Then each $F \in  O_+(4,1)$  maps null lines to null
lines hence preserves the light cone $\mathcal{L}$.  Moreover   it is easy to see that $O_+(4,1)$   acts transitively on $\bb{S}^3$   by  $g. [x] = [gx]$. Here $O_+(4,1)$ is referred to as the group of Moebius transformations of the conformal sphere $\bb{S}^3$. Note that  the subgroup of $O(4,1)$ preserving $P(\mathcal{L}_+)$, is precisely  $O_+(4,1)$.

A smooth map into the conformal sphere $\psi : \Sigma \to \bb{S}^3 \equiv P(\mathcal{L})$ can be viewed as  a
null line subbundle $\Lambda $ of the trivial bundle $  \Sigma  \times  \R^5_1$  via
$\psi(x)= \Lambda_x$, $x \in \Sigma$. 
A (local) lift of $\psi $ is a smooth map $ X: U \to \mathcal{L}$ from an open subset $U \subset \Sigma$, such that the null line spanned by $X(x)$ is $\Lambda_x $ for every $x \in U$.                                                                                
The map   $\psi$ is  called a {\it conformal immersion}  if  every local lift $X$ of $\psi$ is conformal, i.e.          
$\la X_z, X_z \ra =0,  \la X_z, X_{\bar{z}} \ra >0$, for every coordinate $z$. \\

 Let $V: = span\{ X, dX, X_{z \bar{z}} \}$, where $X$ is  a conformal lift of $\psi$.   It is   easily seen  that  $V$ is  in fact independent on the election of a local coordinate $z$ and any  particular conformal lift of $\psi$. So $V$ is can be viewed  as a 
  vector sub-bundle  $V \subset \R^5_1 \times \Sigma$ on which 
  the ambient metric of $\R^5_1$  induces a vector bundle metric of signature $(3,1)$. Each fiber  $V_x$  determines a Moebius invariant $2$-sphere $\bb{S}^2(x) \equiv P(V_x \cap \mathcal{L}) \subset P(\mathcal{L}) \cong \bb{S}^3$.   These spheres altogether  comprise the so-called  {\it mean curvature sphere} or {\it central sphere congruence} of the surface $\psi$~\cite{burstall-pedit-pinkall}.  \\
Respect to  a fixed a local coordinate $z:U \to \C$  there is  a distinguished  local lift $Y: U \to \mathcal{L}_+$ of $\psi$ taking values in the future light cone  such that 
$$  \la Y_z, Y_{\bar{z}} \ra = \frac{1}{2},$$
or equivalently  $|dY|^2 = |dz|^2$ on $U$. It  is called the {\it canonical lift} of the surface $\psi$ and is Moebius invariant. \\
The complementary orthogonal line sub-bundle $V^{\bot}$ is determined by $\Sigma \times \R^5_1 = V \overset{\bot}{\oplus} V^{\bot}$ and the  connection $D$ on $V^{\bot}$ is just  orthogonal projection of the usual derivative in $\R^5_1$: 
$$D_X v = [d_X v]^{\bot}, \quad v \in \Gamma(V^{\bot}), \quad X \in T\Sigma.$$
   
	Let   $N \in \Gamma(V)$ be the unique section satisfying
$$
    \la N,N \ra = \la N, Y_z \ra =\la N, Y_{\bar{z}} \ra = 0, \quad \la Y,N \ra =-1.
$$ 
Thus $V=span\{ Y, Re(Y_z), Im(Y_z), N \}$  and  it is shown in~\cite{burstall-pedit-pinkall} that the  Moebius invariant frame 
 $\{ Y, Y_z, Y_{\bar{z}}, N \} \subset \Gamma( V \otimes \C)$  satisfies  orthogonally relations given by
\begin{equation}
\begin{array}{l}
\la Y, Y \ra = \la N,N \ra =0,\quad   \la N,Y \ra =-1,\\
  \la Y, dY \ra = \la N, dY \ra = \la dN, N \ra=0,  \\
	 \la Y_z, Y_{z} \ra = \la Y_{\bar{z}}, Y_{\bar{z}} \ra =0, \quad  \la Y_z , Y_{\bar{z}} \ra = \frac{1}{2}.\\  
\end{array}
\end{equation}
A direct consequence of the above equations is that $Y_{zz}$ is orthogonal to $Y, Y_z$ and $Y_{\bar{z}}$ and so there is a unique election of a  local complex function $s$ on $\Sigma$ for which  $Y_{zz} + \frac{s}{2} Y $ is a section of the normal bundle $V^{\bot} \otimes \C$ namely,  $\frac{s}{2} = \la Y_{zz}, N \ra$.  Thus  we arrive at the  fundamental equation of Moebius invariant surface geometry:  
\begin{equation}\label{hillseq}
Y_{zz} + \frac{s}{2} Y = \kappa,
\end{equation}
defining uniquely  the complex valued function $s$ and the section $\kappa$ of $V^{\bot} \otimes \C $, respect to the local coordinate $z$. The function $s$ is interpreted as the {\it schwartzian derivative} of the conformal immersion $\psi$ with respect to $z$,
 and  $\kappa$ is identified with  the  {\it  normal valued Hopf differential } of $\psi$, respect to the  coordinate $z$. By construction $s$ and $\kappa$ are Moebius invariants of the immersion $\psi$ with respect to a given coordinate $z$.\\  

In~\cite{burstall-pedit-pinkall} there is an  interpretation of $\kappa$ in terms of  euclidean invariants of the immersion $\psi$ which we briefly describe: There is a unique conformal immersion  $\widehat{\psi}: \Sigma \to \bb{S}^3 \subset \R^4$ satisfying $ [(\widehat{\psi}(x),1   )] = \psi(x)$, $\forall x \in \Sigma$. Thus   $\phi =(\widehat{\psi}(x),1)$ is a lift of $\psi$, which is called  the euclidean lift of $\psi$~\cite{burstall-pedit-pinkall}. Let  $\nu(\widehat{\psi} )$ denote  the normal bundle of  the immersed surface $\widehat{\psi}$. Then there is a bundle isomorphism $\nu(\widehat{\psi} ) \cong V^{\bot}$ given by
\begin{equation}\label{bundleisomorphism}
v \mapsto   \la v, \widehat{H} \ra (\widehat{\psi},1) + (v,0).
\end{equation} 
where $\widehat{H}$ is the mean curvature vector of $\widehat{\psi}$. Under this isomorphism  $\kappa \in \Gamma(V^{\bot} \otimes \C)$  corresponds to a  complex section  $\widehat{\kappa} \in \nu(\widehat{\psi}) \otimes \C$ satisfying  $ \kappa = \la \widehat{\kappa}, \widehat{H} \ra (\widehat{\psi},1) + (\widehat{\kappa},0)$.  
Using~\eqref{hillseq}  it is shown  that  
$$
\widehat{\kappa} \frac{dz^2}{|dz|} = \frac{ \II^{(2,0)} }{ |d \phi|},
$$
 where $\II^{(2,0)}$ is  the $(2, 0)$-part of the normal bundle valued (euclidean) second fundamental form of $\widehat{\psi}$. In this way $\kappa$, up to the isomorphism~\eqref{bundleisomorphism}, is the trace free part of the second fundamental form, i.e., the normal
bundle valued Hopf differential of $\widehat{\psi}$, scaled by the square root of the $\widehat{\psi}$-induced metric. \\

The following structural equations of a  conformal immersion $\psi: \Sigma \to \bb{S}^3$  were   obtained in~\cite{burstall-pedit-pinkall} from  the above orthogonality conditions: 
\begin{equation}\label{structureG}
\begin{array}{ll}
(i) & Y_{zz} = -\frac{s}{2} Y + \kappa, \\
(ii) & Y_{\bar{z}z} = - \la \kappa, \bar{\kappa} \ra Y + \frac{1}{2}N,\\
(iii) & N_z = -2 \la \kappa, \bar{\kappa} \ra Y_z -s Y_{\bar{z}} + 2 D_{\bar{z}} \kappa. \\
\end{array}
\end{equation} 
  The compatibility  among these  are the following equations, 
		\begin{equation}\label{confgausscodazzi}
\begin{array}{ll}
\text{\sf Conformal Gauss:} 
&  \frac{s_{\bar{z}}}{2} = 3 \la \bar{\kappa}_z, \kappa \ra + \la \bar{\kappa},  \kappa_z \ra, \\
\text{\sf Conformal Codazzi:}
 & Im(  \kappa_{\bar{z}\bar{z}} + \frac{\bar{s}}{2}  \kappa )=0.
\end{array}
\end{equation}

Also when the local coordinate changes from $z$  to $w$ the new invariants $s'$ and $\kappa'$ change according to  
\begin{equation}\label{change-s-k}
\begin{array}{ll}
\kappa' = \kappa \left( \frac{\partial z}{\partial w} \right)^{\frac{3}{2}} \left( \frac{\partial \bar{z} }{\partial \bar{w} } \right)^{-\frac{1}{2}},\\
s' = s \left( \frac{\partial z}{\partial w} \right)^2 + S_w(z),
\end{array}
\end{equation} 
where   the usual Schwartzian derivative of a meromorphic function $ g : \Sigma \to \C $ is given by $
  S_z(g) = ( \frac{g''}{g'} )'- \frac{1}{2} ( \frac{g''}{g'} )^2$.
The importance  of the conformal Gauss and Codazzi's equations  is reflected in the following  fundamental theorem of conformal surface theory,
\smallskip

\begin{theorem}~\cite{burstall-pedit-pinkall}\label{fund.theorem.conf.surf.}
Let $\Sigma$ be a Riemann surface and $\psi_j  : \Sigma \to \bb{S}^3$ be conformal immersed surfaces inducing the same Hopf differentials and the same Schwartzians. Then there is a Moebius transformation $T: \bb{S}^3 \to \bb{S}^3$ with $T \psi_1= \psi_2$.\\
Conversely, let $\kappa$ and $s$ be given data on $\Sigma$ transforming according~\eqref{change-s-k}, which also satisfy the conformal Gauss and Codazzi equations~\eqref{confgausscodazzi}. Then there exists a conformal immersion $x: \Sigma \to \bb{S}^3$ with Hopf differential $\kappa$ and Schwartzian $s$.    
\end{theorem}

\begin{remark}\label{kappastory}
 It is proved in~\cite{burstall-pedit-pinkall} that    $\kappa \frac{dz^2}{|dz|}$  is  a globally defined quadratic differential with values in  $L\otimes  \C$,  where $L$ is the real line bundle
$( K \otimes \bar{K})^{−1/2}$ of densities of conformal weight $1$ over $\Sigma$~\cite{calderbank1}.  Then for any  local coordinate system $(U,z)$,   $\kappa$ is can be viewed just as a local complex function on $U \subset \Sigma$ which transforms  according~\eqref{change-s-k}.     
\end{remark}

\begin{remark}
If a conformal immersion $\psi: \Sigma \to \bb{S}^3$ has $\kappa \equiv 0$, then  the image of $\psi$ is contained in a fixed $2$-sphere $\bb{S}^2 \subset \bb{S}^3$, as  follows from~\eqref{structureG}.  
Considering $\psi $ as a conformal map $\psi : \Sigma \to \bb{S}^2 \equiv \bb{CP}^1$, it is shown in~\cite{burstall-pedit-pinkall} that  $s =   ( \frac{\psi''}{\psi'} )'- \frac{1}{2} ( \frac{\psi''}{\psi'} )^2$ which is the usual Schwartzian derivative of $\psi$.    In this case it is shown that  proved that $s$ uniquely determines  $\psi$ up to transformations of  $PSl(2,\C)$, the Moebius transformation group  of $\bb{CP}^1$.  
\end{remark}

The map $\gamma : \Sigma \ni x \mapsto V(x) $ with values in the Grassmannian  $ G_{3,1} (\R^5_1)$  is called {\it the conformal Gauss map} of the immersion $\psi$~\cite{bryant},~\cite{burstall-pedit-pinkall},~\cite{ejiri},~\cite{palmer1}~\cite{xiangma-thesis}.   
   $\gamma$  induces  a positive definite conformal metric   on $\Sigma$ given by $g_{\gamma}= \frac{1}{4} \la d \gamma, d \gamma \ra = | \kappa |^2 |dz|^2$~\cite{xiangma-thesis}.  The Willmore energy of the conformal immersion $\psi$ is defined as the total area of $(\Sigma,  g_{\gamma})$ and is given by
\begin{equation}\label{willmoreG}
W(\psi)=  \frac{i}{2}\int_{\Sigma} | \kappa|^2 dz \wedge d\bar{z},
\end{equation}
which coincides  (up to a constant multiple) with the Willmore energy of the immersion $\psi$~\cite{burstall-pedit-pinkall}. 
A  conformal immersion $\psi : \Sigma \to \bb{S}^3$ is called a {\it Willmore surface} if it extremizes the Willmore energy functional~\eqref{willmoreG}. It is known~\cite{burstall-pedit-pinkall} that $\psi$ is Willmore iif  its conformal invariants $\kappa$ and $s$ (the  Hopf differential and the Schwartzian derivative) satisfy the following stronger version of the conformal Codazzi's equation:
\begin{equation}\label{willmoreEq}
 \kappa_{\bar{z}\bar{z}}+ \frac{\bar{s}}{2} \kappa=0.
\end{equation}

Spacelike  surfaces in $\bb{S}^4_1$ are related to surfaces in $\bb{S}^3$ and $\R^3$ through a double cover of the  conformal Gauss map $\gamma$, which is  the Bryant's Gauss conformal map $Y$: given  an oriented surface $\psi : \Sigma \to \bb{S}^3$ with mean curvature $H$, the conformal Gauss map $Y_{\psi} $ assigns to a point $x \in \Sigma$ 
the oriented sphere $S(x) \subset \bb{S}^3$  of radius $|H(x)|^{-1}$  in  contact with
the surface at $\psi(x)$. Thus  $Y_{\psi}$   takes  values in the manifold of all oriented $2$-spheres (and planes) in $\bb{S}^3$ which is  identified with De Sitter $4$-space $\bb{S}^4_1$~\cite{jeromin}. 
   B. Palmer~\cite{palmer2} observed  that when $Y_{\psi}$ has non-zero mean curvature vector  it  is   marginally trapped. This fact is also implicit in the work of Blaschke~\cite{blaschke}.

\section{ the null gauss map and its conformal invariants   }\label{KonformalinvariantsMTS}
		
Let $f: \Sigma \to \bb{S}^4_1$ be a conformal (spacelike)  immersion  and fix an orientation on the normal bundle $\nu(f)$. Let  $\{ N_1, N_2 \} \subset \Gamma(\nu(f)) $ be a positively oriented  lorentzian  orthonormal frame.  Then for each $p \in \Sigma$  the frame $\{ N_1, N_2 \}$  determines the   null line $span\{ N_1(p)+ N_2(p)\} $. We claim that this null line depends only on $p$ and not  on  $\{ N_1, N_2 \}$.   
 In fact, if    $\{ N'_1, N'_2 \} \subset \Gamma(\nu(f))$ is    another  positively oriented orthonormal   frame then  both frames are related by a gauge,   
\begin{equation*}
\begin{array}{ll}
N'_1 = \cosh (s) N_1 + \sinh (s) N_2,\\
N'_2 = \sinh (s) N_1 + \cosh (s) N_2.\\
\end{array}
\end{equation*}
from these equations it follows that  $ N'_1 + N'_2 = e^s (N_1+ N_2) $,  and so  $N'_1 + N'_2$ and $N_1 + N_2$ generate the same null line. 
Let       
\begin{equation}\label{confGauss}
 G : \Sigma \to \bb{S}^3, \quad   G(x) = [ N_1(x) + N_2(x) ], \,\, x \in \Sigma.
\end{equation}
i.e. $G(x)$ is the null line generated by $N_1(x) + N_2(x)$,  
where $\{ N_1, N_2 \} \subset \Gamma(\nu(f))$ is any positively oriented orthonormal frame. Thus $G$ is  well defined by our previous observation and we call it {\it the null Gauss map of} $f$.  \\

Denote by       $\widehat{G}$  the unique  smooth map from $\Sigma$ to the round euclidean sphere $\bb{S}^3 \subset \R^4$  such that $[(\widehat{G}(x),1) ] =G(x)$ for every $x \in \Sigma$, then    $\phi=:(\widehat{G},1) : \Sigma \to \mathcal{L}_+$ is called     {\it the euclidean lift} of $G$. Thus  $\phi=:(\widehat{G},1)$  takes values in the conic section $
         \mathcal{S}  = \{ x \in \mathcal{L} : \la x ,e_4 \ra =-1 \}$ which   
 inherits from the ambient $\R^5_1$ a positive definite metric of constant curvature $+1$,  and so it is a copy of the round $3$-sphere of radius one. \\

For any  positively oriented orthonormal frame $\{ N_1, N_2 \} \subset \Gamma(\nu(f))$,   $X = N_1+ N_2$ is a  local  lift  of  $G$ with values in $\mathcal{L}_+$. 
 Using  the structure equations~\eqref{structureEqs}   we  see that 
 \begin{equation}\label{derivzgauss}   
   X_z= N_{1,z} + N_{2,z} = (h_2-h_1) f_z + e^{-2u} (\xi_2-\xi_1)  f_{\bar{z}}+ \sigma X. 
   \end{equation}   
 Hence $\la X, f_z \ra = \la X, f_{\bar{z}} \ra =0$. Moreover  since 
 $\la f_z , f_z \ra = \la f_{\bar{z}}, f_{\bar{z}} \ra =0$, $\la f_z, f_{\bar{z}} \ra = e^{2u}$, then       
\begin{equation}\label{derivzgauss2}
\la X_z, X_z \ra = \la N_{1,z} + N_{2,z}, N_{1,z} + N_{2,z} \ra =  (h_2-h_1)(\xi_2-\xi_1)
\end{equation}
Thus if     $f$ is marginally trapped   $\la X_z, X_z \ra =0$.
Let  $Z$ be another  local lift of $G$, then $X= \lambda Z$ for some smooth  non-zero function $\lambda$. Respect to a local coordinate $z$ we compute      
$X_z = \lambda_z Z  + \lambda Z_z $. Since   $\la Z,Z\ra = 0$,  then $0= \la Z, Z_z\ra$, hence   
$     0=  \la X_z, X_z \ra = \lambda^2 \la Z_z, Z_z \ra$, from which $\la Z_z, Z_z \ra =0$  follows. 
On the other hand since    $X_{\bar{z}} = \lambda_{\bar{z}} Z  + \lambda Z_{\bar{z}} $, then from   $\la Z, Z_{\bar{z}} \ra =0$, and~\eqref{derivzgauss}   
we obtain 
\begin{equation}\label{metricGauss}
\lambda^2 \la Z_z, Z_{\bar{z}} \ra =  e^{-2u} | \xi_1-\xi_2|^2 = \la X_z, X_{\bar{z}} \ra.
\end{equation} 
Hence away from the zeros of  $\xi_1-\xi_2$   it follows that $ \la X_z, X_{\bar{z}} \ra >0 $ and  $\la Z_z, Z_{\bar{z}} \ra >0$. In particular if $\xi_1-\xi_2 $ is never zero on $\Sigma$ then  $G :  \Sigma \to \bb{S}^3$ is a conformal immersion.        \\
We call  $q:=(\xi_1-\xi_2 )dz^2$  {\it the Hopf quadratic differential} of the marginally trapped surface $f: \Sigma \to \bb{S}^4_1$.   
The  quadratic Hopf differential  was introduced in~\cite{aledo-galvez-mira} for  marginally trapped surfaces in $\R^4_1$.   
We have proved the  following Lemma: 
\begin{lemma}\label{Gconfimmersion}
Let $f: \Sigma \to \bb{S}^4_1$ be a conformally immersed marginally trapped surface and $q$ its  quadratic Hopf differential.
Then  every (local) lift $Z$ of the null Gauss map $G$ satisfies $\la Z_z, Z_z \ra =0$ and $\la Z_z, Z_{\bar{z}} \ra >0 $  away   the zeros of $q$. In particular  if $q(x) \neq 0, \forall x \in \Sigma$ then   $G: \Sigma \to \bb{S}^3$ is a conformal immersion. 
\end{lemma} 
Since    $(\widehat{G}, 1)$ is a lift of $G$, then  away the zeros of $q$,   $\widehat{G}$ satisfies $\la \widehat{G}_z, \widehat{G}_z \ra =0$ and $\la \widehat{G}_z, \widehat{G}_{\bar{z}} \ra >0$, where $\la.,.\ra $ is the round metric on the sphere $\bb{S}^3$.  Thus if $q$ is never zero    $\widehat{G}$ is  a conformal immersion into the round $3$-sphere. \\

Let    $f: \Sigma \to \bb{S}^4_1$ be  a spacelike immersion then   from   Ricci's equation  $\nu(f)$ is flat if and only if  $Im(\sigma_{\bar{z}})=0$. In this case     
$\sigma_{\bar{z}} - \overline{\sigma_{\bar{z}}}  =  \sigma_{\bar{z}} - \overline{\sigma}_{\,z} = 0$ which shows that the real  one form $\eta := \sigma dz + \overline{\sigma} d \bar{z}$   is closed. Hence     
   there is a locally defined smooth real  function $\beta$ such that $d \beta = \eta$.  One can define  a  new positively oriented orthonormal lorentzian  frame $\{ N'_1, N'_2 \}$  by  
          $$ N'_1 = \cosh(\beta) N_1 + \sinh(\beta)N_2, \quad N'_2 = \sinh(\beta) N_1 + \cosh(\beta)N_2.$$     
Then  it is easy to check that  the new frame $\{ N'_1, N'_2 \}$ has  structure function $\sigma'=0$, so that       $\{ N'_1, N'_2 \}$ is a  $\nabla^{\bot}$-parallel frame which  is unique up to (constant) hyperbolic  rotations  in  $\nu(f)$.   We  keep denoting by $\{ N_1, N_2 \}$ this  new  positively oriented  $\nabla^{\bot}$-parallel orthonormal frame.  
If $f$ is  marginally trapped then   Codazzi's equations~\eqref{gausscodazziricci} reduce to     
\begin{equation}\label{codazziModified}
\xi_{1,\bar{z}}  = \xi_{2,\bar{z}}  = e^{2u} h_{z}, \quad  h= h_1 = h_2,
\end{equation}
which imply   $ (\xi_{1}- \xi_{2})_{\bar{z}}=  e^{2u} (h- h)_z=0$, hence  $q$ is holomorphic. 
Conversely, if $q$ is holomorphic then again by Codazzi's equation we obtain 
  $0= (\xi_1 - \xi_2)_{\bar{z}} = \bar{\sigma} (\xi_1 - \xi_2)$.   
 If  $q$ does not vanish identically then  $\sigma$ must be zero away the isolated zeros of $q$, thus $\sigma\equiv 0$  by continuity.   
  We have  proved the following,
\begin{lemma}\label{holomorphicHopfdiff}
Let  $f: \Sigma \to \bb{S}^4_1$ be a marginally trapped surface. If   $f$  has  flat normal bundle the Hopf   differential $q= (\xi_2- \xi_1) dz^2$ is holomorphic. Conversely, if $q$ is holomorpic and non-identically zero, then $f$ has flat normal bundle.    
\end{lemma}

\begin{remark}\label{gaussmapstationary}
i) If a conformally immersed surface $f: \Sigma \to \bb{S}^4_1$  has zero mean curvature vector   then its normal bundle is not necessarily flat. In this case the Hopf differential $q$ is holomorphic as consequence of Codazzi's equations~\eqref{gausscodazziricci}. \\ 
ii) If  $f: \Sigma \to \bb{S}^4_1$ is marginally trapped with  parallel  mean curvature vector     then $\nu(f)$ is flat~\cite{elghanmi} and so  $q$ is holomorphic by Lemma~\ref{holomorphicHopfdiff}.  \\
iii) From~\eqref{MCV} the $\nabla^{\bot}$-derivative of the mean curvature vector of a marginally trapped surface in a positively oriented normal frame is given by 
 $$\nabla^{\bot}_{\Pz} \h = (h_z+ \sigma h) (N_1+N_2). $$
Thus  $\nabla^{\bot} \h =0$ implies $\nu(f)$ is flat~\cite{elghanmi}, hence  $h$ is constant in a positively oriented $\nabla^{\bot}$-parallel frame $\{ N_1, N_2 \} \subset \Gamma(\nu(f))$.    Conversely if $\nu(f)$ is flat, then $\sigma =0$ for any   $\nabla^{\bot}$-parallel  orthonormal frame $\{ N_1, N_2 \} \subset \Gamma(\nu(f))$. 
\end{remark}

\begin{remark}\label{constantG}
If  $q \equiv 0$, then by~\eqref{metricGauss}   $ N_1 + N_2$ is  a constant null line for  every  oriented lorentzian frame $ \{ N_1 , N_2 \}$, hence the  null Gauss map $G$ is  constant.  Since  $\la f, N_1+N_2 \ra =0$,   the surface $f$ has constant curvature $K=1$ by~\eqref{gausscurvature} and lies in the degenerated hypersurface $M_0 \subset \bb{S}^4_1$,  which is the intersection of the   degenerate $4$-plane $[N_1+N_2]^{\bot}$ in $\R^5_1$ with $\bb{S}^4_1$. For instance this is just the case of any  marginally trapped surface  $f: \bb{S}^2 \to \bb{S}^4_1$ with flat normal bundle. In fact since $q$ is holomorphic on $\bb{S}^2$, it  must vanish. 
 \end{remark}

\subsection{Spacelike isothermic surfaces} \label{isothermicSpaceLike}
The normal valued quadratic Hopf differential of a spacelike immersion $f: \Sigma \to \bb{S}^4_1$ is the $\Gamma(\nu(f) \otimes \C) $-valued two-form 
$$
    \Omega = \xi_1 N_1 dz^2 + \xi_2 N_2 dz^2,
$$
defined in terms of an orthonormal frame $\{N_1, N_2\} \subset   \Gamma(\nu(f))$, 
where $\xi_1, \xi_2$ are the coefficients of $\II( \Pz, \Pz)$, the  $(2,0)$-component   of the second fundamental form of $f$.  
  The spacelike surface  $f: \Sigma \to \bb{S}^4_1$ is called isothermic~\cite{peng-wang} if for each point $x \in \Sigma$ there is a coordinate $z$ for which the normal valued Hopf differential $\Omega$ is real-valued.   Note that from Ricci's equation~\eqref{gausscodazziricciBasic} it follows that every isothermic spacelike immersion in $\bb{S}^4_1$ has flat normal bundle.

\subsection{Non-isotropic spacelike surfaces}\label{nonisotropic}  A     conformally (hence spacelike)  immersed  surface  $f:\Sigma \to \bb{S}^4_1$  is called {\it non-isotropic} if the quartic complex differential $Q=\la f_{zz}, f_{zz} \ra dz^4$ is never zero on $\Sigma$. The quartic complex differential $Q$  was introduced    in~\cite{bryant} in the context  of the conformal Gauss map. In terms of an orthonormal frame $\{N_1, N_2\} \subset \Gamma(\nu(f))$,  $Q=(\xi_1^2 - \xi_2^2)dz^4$, thus   if $f$ is non-isotropic then the Hopf differential $q=(\xi_1-\xi_1) dz^2$ is never zero and  so  the null Gauss map $G:\Sigma \to \bb{S}^3$ is a conformal immersion.   
 The notion of isotropy has an   interpretation  in terms of the {\it curvature hyperbola} which is the image of the unit circle on $T_p \Sigma$ under the second fundamental form of $f$:
$$
      \{ \II_p(X,X) : X \in T_p \Sigma,  \|X\|^2 =1\} \subset T_p^{\bot} \Sigma
$$  
It is shown that    $f$ is non-isotropic if and only if the curvature hyperbola  at each point of $\Sigma$ is non-equilateral~\cite{elghanmi}. 
A conformal non-isotropic spacelike   immersion $f : \Sigma \to \bb{S}^4_1$ with zero mean curvature vector  is also  called harmonic superconformal~\cite{hulett}.  Hence non-isotropic marginally trapped surfaces can be viewed as natural generalizations of harmonic superconformal surfaces. \\

\subsection{ Sphere congruences} \label{conformalinvariants}
Let $f: \Sigma \to \bb{S}^4_1$ be a non-isotropic marginally trapped surface with null Gauss map $G$  and consider the central sphere  congruence of the surface $G: \Sigma \to \bb{S}^3$, given by  the subbundle   $V=span\{ X, dX, X_{z \bar{z}}  \} \subset \Sigma \times \R^5_1$, where $X: \Sigma \to \mathcal{L}_{+}$ is any  local lift of $G$.  Since    $\bb{S}^4_1$ identifies with the manifold of oriented $2$-spheres in $\bb{S}^3$, the immersion $f$ is associated  to    the $2$-sphere congruence   $\Sigma \ni x \mapsto S(x)$,  where $S(x)$ is the $2$-sphere  obtained by projectivization of the intersection of    the Minkowski vector subspace $f^{\bot}(x) \subset \R^5_1$ with the  null cone $\mathcal{L}$: 
$$
S(x) = P( f^{\bot}(x) \cap \mathcal{L}) \subset \bb{S}^3.
$$
Note that the antipodal surface $(-f)$ determines the same sphere congruence $x \mapsto S(x)$.   We say that $S(x)$ is oriented if it is associated to $f$, and opposite oriented if it is associated to $-f$.    
 We claim that $f^{\bot}=V$,  i.e. both sphere congruences coincide. To prove the claim we use the  local  lift of $G$ given by $X:=N_1+N_2 : U \to \mathcal{L}_+$, where $\{N_1, N_2\} \subset \Gamma (\nu(f))$ is a positively oriented orthonormal lorentzian frame.
Thus  $V= span \{ X, Re(X_z), Im(X_z) , X_{z \bar{z}} \}$.  In particular   $\la X,f \ra =0$ since $N_1, N_2$ are normal to $f$.  
On the other hand  from~\eqref{derivzgauss},
\begin{equation}\label{derivzphi1}
         X_z =  e^{-2u}(\xi_2 - \xi_1) f_{\bar{z}} + \sigma X. 
\end{equation}
Hence  $\la f, X_z \ra = \la f, X_{\bar{z}} \ra =0$, or  $\la f , d X \ra=0$. 
Since every lift $ W$ of $G$ is a multiple of $X$ by some function, then  $W$ satisfies $\la f, W \ra =0$ and $\la f , d W \ra=0$. This just says that $G$ is  an envelope of the congruence determined by $f$~\cite{jeromin}. \\
On the other hand taking  $\Bz$ on~\eqref{derivzphi1} and using again~\eqref{structureEqs} yields 
\begin{equation*}
\begin{array}{ll}
X_{z \bar{z}} = e^{-2u} (\xi_2 - \xi_1) ( \bar{\xi}_1 N_1 +  \bar{\xi}_2 N_2)+ 
  \sigma e^{-2u} (\bar{\xi}_2 - \bar{\xi}_1) f_z + ( \sigma_{\bar{z}} + |\sigma|^2) X, \\
\end{array}
\end{equation*}
from which    $\la f ,  X_{z \bar{z}} \ra =0$ follows and so  
       $V \subseteq f^{\bot}$. Thus  $V = f^{\bot}$ since $V$ has rank four. 
        We have proved  the  following 

\begin{proposition}
Let $f: \Sigma \to \bb{S}^4_1$ be a non-isotropic conformal marginally trapped immersion with null Gauss map $G: \Sigma \to \bb{S}^3$.  Then $G$ is an envelope of the spherical congruence determined by $f$. 
Moreover,   the central sphere congruence of  the null Gauss map   $G$  coincides with the spherical congruence determined by $\pm f$.       
\end{proposition}

Recall from Section~\ref{conformalsurfacetheoryS3} that the  correspondence    $ \gamma_G : x \mapsto V(x) $ defines the {\it conformal Gauss map} of the surface $G: \Sigma \to \bb{S}^3$. Since $V= f^{\bot}$, then   $\gamma_G$   takes values in $G_{3,1} (\R^5_1)$  the Grassmannian of all subspaces of $\R^5_1$ with signature $(+++-)$. Since $V$ and $V^{\bot}=\R f$ determine each other then either of them can be used to define the conformal Gauss map of $G$. Thus for each $x\in \Sigma$, $\gamma_G(x) = \R f(x)$ belongs to the   
 manifold of all spacelike lines through the origin of $\R^5_1$ which identifies also with  $G_{3,1} (\R^5_1) $.  Note that the projection     $\bb{S}^4_1 \to G_{3,1} (\R^5_1)$ given by $P: p \mapsto \R p$ is a lorentzian double cover. Intersecting  the spacelike line $\gamma_G(x) = \R f(x)$ with  $\bb{S}^4_1$ we obtain $ \{ +f(x), -f(x) \} \subset \bb{S}^4_1$ which is just   the fiber of $P$ over $G(x)\in \bb{S}^3$.  Thus the surface $f$ and its antipodal $-f$ have the same null Gauss map $G$. Thus  the null Gauss map $G$ can be considered as a pseudo-inverse of the conformal Gauss map $\gamma_G$.\\

\subsection{An equation  relating $\kappa, s$ and $\delta$} Let $Y$ be the canonical lift of $G$ respect to a local coordinate $z$.  Then   there is  a  non-zero function $\tau$ such that $X = \tau Y$. Using~\eqref{derivzgauss}, we compute
\begin{equation}\label{derivzphi2}
         \tau_z Y + \tau Y_z = X_z = e^{-2u}(\xi_2 - \xi_1) f_{\bar{z}}+ \sigma X. 
\end{equation}
Hence $\la X_z, X_{\bar{z}} \ra = \frac{\tau^2}{2} = \tau^2 \la Y_z, Y_{\bar{z}} \ra = e^{-2u} |\xi_2 - \xi_1|^2$, so that     
\begin{equation}\label{functiont}
\tau =    \sqrt{2} e^{-u} |\xi_2 - \xi_1|.
\end{equation}
 Hence we   obtain the  canonical lift of $G$ in terms of $X = N_1 +N_2$:
$$
           Y = \frac{e^{u}}{\sqrt{2}|\xi_2 - \xi_1|} (N_1+N_2).
$$  
A routine computation using the structure equations of $f$ shows  that $Y$ is in fact independent   on any particular choice of a positively oriented lorentzian frame $\{ N_1, N_2\}$. 
On the other hand  
\begin{equation*}
\begin{array}{l}
 \tau_{zz} Y + 2\tau_z Y_z + \tau Y_{zz} = X_{zz} =\\
\\
 (e^{-2u}(\xi_2 - \xi_1))_z f_{\bar{z}} + e^{-2u}(\xi_2 - \xi_1) (-e^{2u} f + e^{2u} h X) + \\
\\
\sigma_z X + \sigma \{ e^{-2u}(\xi_2 - \xi_1) f_{\bar{z}} + \sigma X  \}.
\end{array}
\end{equation*}
 
Adding and substracting $\tau \frac{s}{2}Y $ we obtain
\begin{equation}\label{decomp1}
\begin{array}{c}
 (\tau_{zz} -\tau \frac{s}{2})Y +  2 \tau_z Y_z + \tau (Y_{zz}+\frac{s}{2} Y)  =\\
\\
 (e^{-2u}(\xi_2 - \xi_1))_z f_{\bar{z}} + e^{-2u}(\xi_2 - \xi_1) (-e^{2u} f + e^{2u} h X)+\\
\\
\sigma_z X + \sigma \{ e^{-2u}(\xi_2 - \xi_1) f_{\bar{z}} + \sigma X  \}.
\end{array}
\end{equation}

Comparing  the $V^{\bot}$ components in this identity we obtain the equality     
\begin{equation}\label{kappaG}
 ( \xi_1 - \xi_2) f =  \tau  (Y_{zz}+\frac{s}{2} Y) = \tau \kappa, 
\end{equation}    
 Inserting the function $\tau$ of ~\eqref{functiont}  we obtain 
a  formula for the normal valued Hopf  differential of $G$ which makes sense only if the Hopf quadratic deifferential of $f$ is non-zero:
\begin{equation}\label{ExpressionKappa}
     \kappa = \frac{(\xi_1- \xi_2) e^{u}}{\sqrt{2} | \xi_1 - \xi_2|}f. 
\end{equation}
Using the polar form   
$(\xi_1- \xi_2)=  | \xi_1 - \xi_2|  e^{i \theta} $  the above expression becomes  $\kappa =  \frac{e^{u + i \theta}}{\sqrt{2}} f $ and so by  Remark~\ref{kappastory} we  identify     
\begin{equation}\label{ExpressionKappa2}
\kappa \equiv  \frac{e^{u + i \theta}}{\sqrt{2}}, \quad \text{where} \,\, \frac{(\xi_1- \xi_2)}{ | \xi_1 - \xi_2|}=    e^{i \theta}.  
\end{equation}
In particular we recover the conformal parameter from $\kappa$ above:
\begin{equation}\label{kappaVsu}
  e^{2u} =  2 \la \kappa , \overline{\kappa} \ra.
\end{equation}

In~\cite{burstall-pedit-pinkall} it is shown that  any  section   $v \in \Gamma (V \otimes \C)$ can be decomposed  as follows: 
\begin{equation}\label{sectioninV}
   v = -\la v, N \ra Y - \la v, Y \ra N + 2 \la v, Y_{\bar{z}} \ra Y_z + 2 \la v, Y_{z} \ra Y_{\bar{z}}.
\end{equation}
We use this  formula  to expand the particular section   $f_z \in \Gamma(V \otimes \C)$. 
Since $\tau Y = N_1+ N_2=X$, it follows  $\la f_z, Y \ra=0$.  
Also from  $0= \la f, Y_z \ra_z = \la f_z, Y_z \ra + \la f, Y_{zz} \ra$, equation ~\eqref{structureG}-$(i)$, and  $\la f, Y \ra =0$, we compute   
$$
   \la f_z, Y_z \ra = - \la f, Y_{zz} \ra = -\la f, -\frac{s}{2} Y + \kappa \ra= - \la f, \kappa \ra = - \frac{e^{u + i \theta}}{\sqrt{2}}.  
$$  

On the other hand since   $0 = \la f, Y_{\bar{z}} \ra_z = \la f_z, Y_{\bar{z}} \ra + \la f, Y_{z \bar{z}} \ra$, then 
$$
\la f_z, Y_{\bar{z}} \ra =  -\la f, Y_{z \bar{z}} \ra = | \kappa |^2 \la Y,f \ra - \frac{1}{2} \la N,f \ra =0.
$$ 
Also  $\la f, N \ra =0$, implies   $\la f_z, N \ra + \la f, N_z \ra =0$. 
Hence  $\la f_z, N \ra = - \la f, N_z \ra = -2 \la f, D_{\bar{z}} \kappa \ra $. Since $D_{\bar{z}} \kappa = (u+i \theta)_{\bar{z}} \kappa$, then    
$$
 \la f_z, N \ra  =   -\sqrt{2} (u+i \theta)_{\bar{z}}  e^{u + i \theta}.  
$$

From  these equations and using~\eqref{sectioninV} with $v = f_z$, we obtain   

\begin{equation} \label{fz2}
f_z = \sqrt{2} e^{u + i \theta} \{ (u+i \theta)_{\bar{z}} Y  -  Y_{\bar{z}} \}.
\end{equation}

Therefore,  
\begin{equation}\label{formulaf1}
f_{z\bar{z}} =  \sqrt{2} e^{u + i \theta} \{  ( (u+i \theta)_{\bar{z}}  )^2 + (u+i \theta)_{\bar{z} \bar{z}} +\frac{\bar{s}}{2} \}  Y -  \sqrt{2} e^{u + i \theta} \bar{\kappa}.
\end{equation}

On the other hand using  the structure equations of the immersion $f$ and  $X=N_1+ N_2 = \tau Y$, we obtain     
\begin{equation}\label{formulaf2}
f_{z\bar{z}} = - e^{2u} f + e^{2u} h X= - e^{2u} f + e^{2u} h \tau Y.   
\end{equation}  

Note that    $\sqrt{2} e^{u + i \theta} \bar{\kappa}= e^{2u}f$, so that   equating~\eqref{formulaf1} and~\eqref{formulaf2} gives  
$$
 e^{2u}  h \tau = \sqrt{2}  e^{u + i \theta} \{ ( (u+i \theta)_{\bar{z}}  )^2 + (u+i \theta)_{\bar{z} \bar{z}} + \frac{\bar{s}}{2} \}. 
$$  

Inserting the function $\tau$ given by~\eqref{functiont} in this expression we obtain the following formula:
\begin{equation}\label{schwartzian1}
 h | \xi_2 - \xi_1|e^{-i \theta} =   ((u+i \theta)_{\bar{z}})^2 + (u+i \theta)_{\bar{z} \bar{z}} + \frac{\bar{s}}{2},   
\end{equation}
or conjugating  both sides,
\begin{equation}\label{schwartzian2}
 h (\xi_1- \xi_2)  =   ((u-i \theta)_{z})^2 + (u-i \theta)_{zz} + \frac{s}{2}.
\end{equation}
Now recall   the connection $D$ on the normal bundle $V^{\bot}$. Any section $v \in \Gamma(V^{\bot})$ can be written as  $v=b f$ for some smooth function $b$. Thus $d_X ( b f) = d_X b f + b d_X f$. Condition $df \bot f$ implies $D_X f=0$, hence  
\begin{equation}\label{thenormalconnectionD}
D_X (v)= (d_X b) f.
\end{equation}
Thus we may identify $D_X (v) \equiv  d_X b$. 
 Since   $\kappa \equiv \frac{e^{u + i \theta}}{\sqrt{2}} $, we compute    
$$
  D_{\bar{z}} D_{\bar{z}} \kappa  = \kappa_{\bar{z}\bar{z}} = \left( (u+i\theta)_{\bar{z}}^2 + (u+i\theta)_{\bar{z}\bar{z}} \right) \kappa.
$$
On the other hand since         $\overline{ h(\xi_1- \xi_2)} \kappa = \frac{e^{u}}{\sqrt{2}} h | \xi_2 - \xi_1 | $, then    $ \overline{ h(\xi_1- \xi_2)}  \kappa$ is  real valued and so     equation~\eqref{schwartzian1} becomes 
\begin{equation}\label{Conf.fund.equation}
                 \kappa_{\bar{z} \bar{z}}  + \frac{\bar{s}}{2} \kappa  =     Re \left( \overline{h(\xi_1- \xi_2) }\, \kappa \right).
\end{equation}

Equation  \eqref{Conf.fund.equation}  relates  the quadratic  differential  $h(\xi_1- \xi_2)dz^2$ of a marginally trapped surface  $f:\Sigma \to \bb{S}^4_1$ and   the conformal invariants $\kappa, s$ of its null Gauss map $G$. 
Since   the quadratic differential $\delta:=h(\xi_1- \xi_2) dz^2$   plays a key  role in~\eqref{Conf.fund.equation}, we call   it the $\delta$-differential of the marginally trapped surface $f$.

\begin{remark}
 Equation~\eqref{Conf.fund.equation} implies the conformal Codazzi equation $Im(\kappa_{\bar{z} \bar{z}}  + \frac{\bar{s}}{2} \kappa  )=0$. The conformal Gauss equation~\eqref{confgausscodazzi} may be recovered from~\eqref{schwartzian2} by a long   calculation using   Gauss, Codazzi and Ricci's equations~\eqref{gausscodazziricci}.\\
\end{remark}

\subsection{Congruence} A    basic question is  to what extent  a  marginally trapped surface is  determined   by the conformal invariants  of its null Gauss map. We prove,
   
\begin{theorem}\label{KSdetermineDelta}
Let $f,f' : \Sigma \to  \bb{S}^4_1$ be  non-isotropic  marginally trapped surfaces with null Gauss maps $G,G'$. If $\kappa = \kappa', s = s'$   then there is an isometry $\Phi $ of $\bb{S}^4_1$ such that  $\Phi f = f'$. As a  consequence of this $\delta = \delta'$. 
  
\end{theorem}
\noindent {\bf Proof.} 
By Theorem~\ref{fund.theorem.conf.surf.} there is a Moebius transformation $T \in O_+(4,1)$ of $\bb{S}^3 $ such that $TG=G'$. 
 Recall that the  Moebius group  $O_+(4,1)$ acts on $\bb{S}^3$  by $T([x]) = [Tx]$,  $\forall x \in \mathcal{L}$. 
 Let   $Y$ be   the canonical lift of $G$ respect to to a holomorphic coordinate $z$, then $Y'= TY$ is  the canonical lift of $G'$ respect  to $z$. Since $V = span\{ Y, Re(Y_z), Im(Y_z) , Y_{z \bar{z}} \}$, it follows that  $TV = V'$ and so  $T V^{\bot} = V'^{\bot}$.   This last equality implies  $Tf = \pm f'$ where the sign ambiguity  reflects the fact  that the sphere congruences determined by   $f$ and $f'$   are (modulo Moebius transformations) equal up to  orientation. Defining  $\Phi = T$, if $Tf = f'$ and  $\Phi =  -T$, if $Tf=-f'$, then  $\Phi$ is an isometry of $\bb{S}^4_1$ satisfying $\Phi  f = f'$.
In particular if $T$ is the identity, then $G=G'$ and so   $V^{\bot} = \R f = \R f'$, which implies    $f' = \pm f$.\\
Let $\{ N_1, N_2 \} \subset \Gamma(\nu(f))$ be a positively oriented orthonormal frame, then $\{ \Phi N_1, \Phi N_2 \} \subset \Gamma(\nu(f'))$ is an orthonormal frame. We  can choose an orientation on $\nu(f')$  so that $\{ \Phi N_1, \Phi N_2 \} \subset \Gamma(\nu(f'))$ is a positively oriented normal frame along $f'$.   
Since  $\h = h(N_1 + N_2)$ is the mean curvature vector of $f$, then $\Phi \h = h (\Phi N_1 + \Phi N_2)$ is the mean curvature vector of $f'$. Also since $\II( \Pz, \Pz) = \xi_1 N_1 + \xi_2 N_2$, then $\II'( \Pz, \Pz) = \xi_1 \Phi N_1 + \xi_2 \Phi N_2$ and so  $\delta' = h(\xi_1 - \xi_2)dz^2 = \delta$.  
\hfill $\square$\\

\noindent As a   partial converse of the previous Theorem we obtain   the following  

\begin{lemma}\label{resultoncongruence}
Let $f,f' : \Sigma \to  \bb{S}^4_1$ be  non-isotropic  marginally trapped surfaces which induce the same conformal metric. If either \\
 i)  $f,f'$  are both non-stationary  and  $\delta = \delta'$,  or \\
ii) $f,f'$ are both stationary with  $q =q'$, \\
then  there is an isometry $\Phi$ of $\bb{S}^4_1$ such that $\Phi \circ f = f'$. 
\end{lemma}

\noindent {\bf  Proof:} 
Assume first that $f,f'$ are both non-stationary with   $\delta = \delta'$ i.e.  $h(\xi_1 - \xi_2)dz^2 = h'(\xi'_1 - \xi'_2) dz^2$, hence $h(\xi_1 - \xi_2) = h'(\xi'_1 - \xi'_2)$.  Since $h,h'$ are  real and non-zero, we may assume they are both positive  (if say $h<0$,  we can  replace  $f$ by its antipodal $-f$ which has mean curvature function $-h>0$). Since by hypothesis the Hopf differentials $q,q'$ are never zero, we use the  polar form  $ \xi_1 - \xi_2 = |\xi_1 - \xi_2|e^{i\theta}$ and $ \xi'_1 - \xi'_2 = |\xi'_1 - \xi'_2|e^{i\theta'}$. Hence the equality $\delta = \delta'$ implies   
$$
 h |\xi_1 - \xi_2| e^{i\theta } =h' |\xi'_1 - \xi'_2| e^{i \theta' }.  
$$
It follows that     $\theta -\theta' = 2 k\pi$ with integer $k$. Since  by hypothesis $f$ and $f'$ induce the same conformal metric, we have $u=u'$ and so~\eqref{ExpressionKappa2} implies $\kappa= \kappa'$. 
On the other hand from  $\delta = \delta'$ and~\eqref{Conf.fund.equation} it follows  that  $s=s'$.
Thus       $G, G'$ have the same conformal invariants $\kappa$ and $s$, hence i) follows by  applying the preceding Theorem. \\
If now $f,f'$ are both stationary with $q=q'$, then  $ |\xi_1 - \xi_2| e^{i\theta } = |\xi'_1 - \xi'_2| e^{i \theta' }$, and so $\theta - \theta'$ is an integer multiple of $2 \pi$. Thus since  $u =u'$ by hypothesis,~\eqref{ExpressionKappa2} implies $\kappa = \kappa'$.  Since $f,f'$ are both stationary, then $ \delta = \delta' =0$. Thus from~\eqref{Conf.fund.equation}, we conclude that $s = s'$, and so  $G,G'$ have the same conformal invariants. 
   \hfill $\square$\\

A conformal immersed surface  $\psi : \Sigma \to \bb{S}^3$ is called {\it constrained Willmore} if it extremizes the Willmore energy functional with respect to variations through conformal immersions~\cite{burstall-pedit-pinkall}. It has been proved in~\cite{bohle} that  $\psi$ is constrained Willmore if and only if its  conformal invariants $\kappa, s $  satisfy 
\begin{equation}\label{constrWillmoreEq}
  \kappa_{\bar{z}\bar{z}} + \frac{\bar{s}}{2} \kappa = Re( \bar{\eta} \kappa),
\end{equation}
for some holomorphic quadratic differential $\eta dz^2$ on $\Sigma$. Equations~\eqref{constrWillmoreEq} and~\eqref{Conf.fund.equation} are  related.  In fact, we have seen before that for an immersed non-isotropic marginally trapped surface $f: \Sigma \to \bb{S}^4_1$ the quantity $\overline{h(\xi_1- \xi_2)}\kappa$ is real, so that we ask under what conditions is $\delta = h(\xi_1-\xi_2) dz^2$ holomorphic.  

\begin{lemma}\label{delta-QholomoprphicNablaH}
Let $f: \Sigma \to \bb{S}^4_1$ be a non-isotropic conformally  immersed marginally trapped surface with non-zero mean curvature vector. Then the following affirmations are equivalent:\\
i) The quartic complex differential $Q = \la f_{zz}, f_{zz} \ra dz^4$ is holomorphic,\\
ii) The quadratic complex differential $\delta = h(\xi_1 - \xi_2) dz^2$ is holomorphic,\\
iii)  $f$ has parallel mean curvature vector.
\end{lemma}
\noindent {\bf Proof.} Let $\{ N_1, N_2 \} \subset \Gamma(\nu(f))$ be a positively oriented orthonormal frame, then the quartic differential  becomes $Q= (\xi_1^2-\xi_2^2) dz^4$, where $\xi_1 = \la f_{zz}, N_1\ra, \xi_2 = -\la f_{zz}, N_2 \ra$ and $\h = h(N_1+N_2)$. 
Since $f$ is marginally trapped  Codazzi's equations~\eqref{gausscodazziricci} reduce to 
$$
e^{-2u} ( \xi_{1 \bar{z}} + \bar{\sigma} \xi_2) =e^{-2u} ( \xi_{2 \bar{z}} + \bar{\sigma} \xi_1) =  h_z + \sigma h.
$$
Using  these quations we compute
\begin{equation*}
\begin{array}{l}
 (\xi_1^2-\xi_2^2)_{\bar{z}} = 2 \xi_1 \Bz \xi_1-2 \xi_2 \Bz \xi_2 = \\
2 \xi_1 ( e^{2u} (h_z + \sigma h)-\xi_2 \bar{\sigma} ) -2 \xi_2 ( e^{2u} (h_z + \sigma h)-\xi_1 \bar{\sigma} )=\\
 2 e^{2u}  (h_z + \sigma h) ( \xi_1- \xi_2).
\end{array}
\end{equation*}
Since $f$ is non-isotropic  $q$ is never zero, so  $Q$ is holomorphic if and only if $h_z + \sigma h=0$, which is just the parallel mean curvature equation~\eqref{normalderivativeH}. This proves i) $\Leftrightarrow $ iii).\\
Again  from Codazzi's equation we get     
   $ (\xi_1 - \xi_2)_{\bar{z}} = \bar{\sigma} (\xi_1 - \xi_2)$, which implies  $(h(\xi_1 - \xi_2))_{\bar{z}} = (h_{\bar{z}}+ \bar{\sigma}h)(\xi_1 - \xi_2)$. Hence    
$(\xi_1^2-\xi_2^2)_{\bar{z}} =  2 e^{2u} (h(\xi_1 - \xi_2))_{\bar{z}}$, thus    		  $\delta$   is holomorphic if and only if $Q$ is holomorphic, and so i) $\Leftrightarrow$ ii).
\hfill $\square$\\

Note  for instance  that there is no non-isotropic spacelike immersion $f:\bb{S}^2 \to \bb{S}^4_1$ with parallel non-zero mean curvature vector.  
Isotropic marginally trapped surfaces  in $\R^4_1$, and $ \bb{S}^4_1$ have been considered  in~\cite{cabrerizo-et-al}.   \\ 

As a first consequence of equation~\eqref{Conf.fund.equation} we deduce that a non-isotropic conformal marginally trapped immersion  $f: \Sigma \to \bb{S}^4_1$  has zero mean curvature vector if and only if its null Gauss map $G: \Sigma \to \bb{S}^3$ is a Willmore surface. For, $\h=0$ if 
and only if 
   $\delta \equiv 0$ by~\eqref{MCV}  if and only if~\eqref{Conf.fund.equation} becomes $\kappa_{\bar{z} \bar{z}}  + \frac{\bar{s}}{2} \kappa =0$ which is just the condition for    $G: \Sigma \to \bb{S}^3$ being  a    Willmore surface. 
We obtain  also the following result as consequence of~\eqref{Conf.fund.equation}:

\begin{theorem}\label{GaussConstWillmore}
Let $f : \Sigma \to \bb{S}^4_1$ be a non-isotropic conformal marginally trapped immersion with null Gauss map $G$ and   mean curvature vector $\h \neq 0$.  Then   $\nabla^{\bot} \h =0$ if and only if    $G: \Sigma \to \bb{S}^3$ is  a  constrained Willmore surface.\\  
\end{theorem} 
										
\noindent {\bf Proof:} 
The conformal invariants $\kappa,s$ and the $\delta$-differential of $f$ satisfy  equation~\eqref{Conf.fund.equation} in which $\overline{h(\xi_1-\xi_2)}\kappa$ is real valued. If $f$ has non-zero parallel mean curvature vector then $\delta=h(\xi_1- \xi_2)dz^2$ is holomorphic by Lemma~\ref{delta-QholomoprphicNablaH}. This precisely  says that  $G: \Sigma \to \bb{S}^3_1$    is  constrained Willmore.\\
Conversely if the null Gauss map $G: \Sigma \to \bb{S}^3_1$ of $f$  is a  constrained Willmore surface then  $
\kappa_{\bar{z} \bar{z}}  + \frac{\bar{s}}{2} \kappa  =   Re \left( \overline{\eta }\, \kappa \right)$,   for some holomorphic quadratic differential $\eta dz^2$.  But    $\kappa, s$ uniquely determine   the $\delta$-differential of $f$ by Theorem~\ref{KSdetermineDelta}, so that    $\delta = \eta dz^2$. Therefore   $\delta$ is holomorphic which implies that $f$ has parallel mean curvature vector  by Lemma~\ref{delta-QholomoprphicNablaH}.   \hfill $\square$\\

\section{one-parameter deformations  and associated families}\label{deformationsMTS}
In classical  minimal surface theory an  interesting problem is to determine whether a given minimal surface can be deformed in a nontrivial way. The oldest known example is the deformation of the catenoid into the helicoid~\cite{spivak}. We consider   here two different  non-trivial one-parameter  isometric deformations of   marginally trapped surfaces in $\bb{S}^4_1$. Throughout we only consider non-isotropic surfaces. 

\subsection{The $\bb{S}^1$-deformation  family of marginally trapped surfaces with non-zero parallel mean curvature} 

In submanifold theory the harmonicity of the   Gauss map characterizes  submanifolds of spaceforms with parallel mean curvature vector. This is referred to as the Ruh-Vilms property after the well known paper~\cite{ruh-vilms}. We obtain  here  integrable one-parameter deformations  of  marginally trapped surfaces in $\bb{S}^4_1$ with non-zero parallel mean curvature determined by the spectral symmetry of the harmonic map equation of the $\partial$-transform  of such surfaces.\\

\noindent  Given a conformal immersion $f: \Sigma \to \bb{S}^4_1$  we  consider the map          
$$
\phi: \Sigma \to \bb{CP}^4, \quad x \mapsto [f_z(x)],
$$  
 where $[f_z]  \subset \C^5_1$ is the  spacelike isotropic complex line generated by $f_z$.    Thus  $\phi$ is well defined since it it independent on the local coordinate $z$ and is called  the $\partial$-transform  of the  surface $f$~\cite{wood}. Since $f$ is spacelike  $[f_z(x)]$ is a spacelike complex line hence it is a point in  $\bb{CP}^4_1$,    the open submanifold of  $ \bb{CP}^4 $ consisting of all spacelike complex lines through the origin of $\C^5_1$. Moreover since $f$ is conformal      $\phi$ factors through the  manifold of isotropic spacelike complex lines in $\bb{CP}^4_1$ which  is  the complex quadric defined by
\begin{equation}\label{quadricQ}
Q=\{[z] \in \bb{CP}^4_1 : z^2_0+z^2_1+z^2_2+z^2_3-z^2_4=0\}.
\end{equation}
 Note that since $Q$ is a complex submanifold,  it is   totally geodesic in $\bb{CP}^4_1$.  
Denote by   $L \to \bb{CP}^4_1$   the tautological line bundle whose fiber over a point $l \in \bb{CP}^4_1$ is the line $l$ itself and consider the complex line sub bundle $\ell:=\phi^*(L) \subset \Sigma \times \C^5_1$. Denote  by $\ell^{\bot}$ the complementary orthogonal line sub bundle so that $\Sigma \times \C^5_1 =  \ell \oplus \ell^{\bot}$.  Any  section $\mu $ of the trivial bundle $\Sigma \times \C^5_1$,  decomposes uniquely as $\mu = \mu_1 + \mu_2$ with $\mu_1 \in \Gamma(\ell)$ and $\mu_2 \in \Gamma(\ell^{\bot})$. The projection maps are  defined by $\pi_{\ell} \mu = \mu_1$ and $\pi_{\ell^{\bot}} \mu = \mu_2$.   
Since $Q$ is totally geodesic in $\bb{CP}^4_1$, the map $\phi$ is harmonic as a map into $Q$ if and only if it is harmonic as a map into  $\bb{CP}^4_1$. 
Consider on $\ell$ and $\ell^{\bot}$ the Koszul-Malgrange complex structure~\cite{eells-lemaire}:  a section $s \in \Gamma(\ell)$ (resp. $s \in \Gamma(\ell^{\bot})$) is holomorphic if and only if $\pi_{\ell}(s_{\bar{z}} ) =0$, (resp. $\pi_{\ell^{\bot}}(s_{\bar{z}} ) =0$).
It is known  that    $\phi: \Sigma \to \bb{CP}^4_1$ is harmonic if and only if the map 
$$d \phi (\Pz): \ell \to \ell^{\bot}, \quad  d \phi (\Pz)\mu = \pi^{\bot}_{\ell} (\Pz \mu) $$
  is holomorphic, i.e. it sends holomorphic sections of $\ell$, to holomorphic sections of $\ell^{\bot}$~\cite{burstall-rawnsley},~\cite{eells-lemaire}.  
Recall from~\eqref{structureEqs2} the structure equations of a conformal immersion    $f : \Sigma \to \bb{S}^4_1$ and the corresponding Gauss Codazzi and Ricci's equations~\eqref{gausscodazziricciBasic}.
The second structure equation    $ f_{\bar{z} z} = -e^{2u} f + e^{2u}\h$  implies     $\pi_{\ell}(\Bz f_{ z}) =0$, which says   that $f_z$ is a holomorphic section of $\ell$. In particular    every holomorphic section of $\ell$ is of the form $\zeta  f_z$, where  $\zeta$ a complex holomorphic function. Thus    $\phi$ is harmonic if and only  if $\mu:=\pi_{{\ell}^{\bot}} (f_{zz}) $ is a holomorphic section of $\ell^{\bot} $. From the second structure equation it follows that    $\mu =  \xi_1 N_1 + \xi_2 N_2 = \II ( \Pz, \Pz) \in \Gamma(\ell^{\bot})$.   From~\eqref{structureEqs2} and Codazzi's equations~\eqref{gausscodazziricciBasic}  we compute  
\begin{equation*}
\begin{array}{l}
\pi_{{\ell}^{\bot}} (\Bz \mu) = \pi_{{\ell}^{\bot}} ( \Bz \xi_1  N_1 + \Bz \xi_2  N_2  + \xi_1  \Bz N_1 +  \xi_2  \Bz N_2) =\\
\pi_{{\ell}^{\bot}} \left (  \Bz \xi_1 N_1 +\Bz \xi_2 N_2  +  \xi_1 (-h_1 f_{\bar{z}} - e^{-2u} \bar{\xi}_1 f_z+ \bar{\sigma} N_2) +  \xi_2 (h_2 f_{\bar{z}} +  e^{-2u} \bar{\xi}_2 f_z+ \bar{\sigma} N_1) \right) = \\
 e^{-2u} ( (\Pz h_1 +\sigma h_2)N_1 + (\Pz h_2 +\sigma h_1) N_2 ) = e^{-2u} \nabla^{\bot}_{\Pz} \h.\\ 
\end{array}
\end{equation*}
Hence  the section $\mu$  is holomorphic if and only if  $f$ has parallel mean curvature vector field and this in turn is equivalent to the harmonicity of     $\phi$. We summarize the above discussion in the following Lemma  which is a manifestation of the characterization due to Ruh-Vilms of submanifods with parallel mean curvature vector in $\R^n$ and $S^n$~\cite{ruh-vilms}: 

\begin{lemma}\label{deltaPhi}
Let $f: \Sigma \to \bb{S}^4_1$ be a conformally immersed surface. Then the  Gauss transform $\phi : \Sigma \to Q \subset \bb{CP}^4_1$ of $f$   is harmonic if and only if the  surface $f$ has parallel mean curvature vector. \\
\end{lemma}  

\begin{remark}\label{substantial}
If  $f$ is a non-isotropic marginally trapped with non-zero parallel mean curvature vector then      the codimension  of $f$ cannot be reduced.   This is fails to be  true if the mean curvature vector of $f$  has non-zero squared norm~\cite{chen-vanderveken}.  \\
\end{remark}

The Lie group   $SO_+(4,1)$ acts transitively by isometries on $Q$, where we consider a multiple of  the Killing metric on $Q$. Fixing for instance   the base point $o:= [e_1-ie_2] \in Q$, then $Q$ is diffeomorphic to the  symmetric quotient  $SO_+(4,1)/H$,  where $H$ is the stabilizer of the base point $o \in Q$. 
 Consider the involutive isomorphism  $\tau$ of $SO_+(4,1)$  given by  $\tau(F) = EFE$, where  $ E:= diag(1,-1,-1, 1,1) \in SO_+(4,1)$.  Then  the connected component $Fix(\tau)_0$ of the subgroup  of fixed points of $\tau$  coincides with $H$ which is isomorphic to $SO(2) \times SO(3,1)$.  
The  $(\pm 1)$-eigenspaces   of $d \tau_e$ are given respectively by 
\begin{equation}\label{eigenspacesTAU}
\goth{m}: = \{   \left( \begin{smallmatrix}
      0  & a & b & 0 & 0 \\
      -a & 0 & 0 & c & d \\
      -b & 0 & 0 & e & k\\
      0  & -c&-e & 0 & 0 \\
      0  &   d&k & 0 & 0 \\
      \end{smallmatrix} \right): a,b,c,d,e,k \in \R \}, \quad     
 \goth{h}: = \{   \left(  \begin{smallmatrix}
      0  & \,0  & \,0 & \,m & n \\
      0  & \, 0 & s   & \,0 & 0 \\
      0  & -s\,   & \,0 & \,0 & 0\\
      -m\,  & \,0  & \,0 & \,0 & t \\
      n  & \, 0 & \,0 & \,t & 0 \\
      \end{smallmatrix} \right): s, t,  m,n  \in \R\},    
\end{equation}
which  satisfy  
$
[ \goth{h}, \goth{m}] \subseteq \goth{m}, \quad [ \goth{h}, \goth{h}] \subseteq \goth{h}, \quad [ \goth{m}, \goth{m}] \subseteq \goth{h}$.\\  

We briefly review the loop group formulation of the harmonic map equation for maps into a symmetric space.    Let $G$ be a connected semisimple (compact or non-compact)  Lie group and  assume  that $G$ is a matrix group.  Let $G/H$ be an inner symmetric space with  involution $\tau : G \to G$ satisfying $ (G_{\tau})_0 \subseteq H \subseteq G_{\tau}$, then $G/H$ has a $G$-invariant  non-degenerate symmetric bilinear form~\cite{helgason}. Let $\goth{g}=Lie(G), \goth{h}=Lie(H)$ be the Lie algebras of $G$ and $H$ respectively. The involution $\tau$ induces a decomposition $\goth{g}= \goth{h} \oplus \goth{m}$ into eigenspaces of $d \tau_e$ such that $\goth{h} =\{X: d \tau_e(X) = -X\}$ and $\goth{m}=\{X: d \tau_e(X) = X\}$. It follows that these eigenspaces satisfy $[ \goth{h}, \goth{h} ] \subset \goth{h}, [ \goth{h}, \goth{m} ] \subset \goth{m}, [ \goth{m}, \goth{m} ] \subset \goth{h} $. \\

Let   $\psi : \Sigma \to G/H$ be a smooth map and  $F: U \to G$  a  frame of $\psi$ on $U$, where $U \subset \Sigma$ is a simply connected open subset (if $\Sigma$ is simply connected then there is always a global frame $F: \Sigma \to G$)  . Let $\alpha:= F^{-1}dF$ be the Maurer-Cartan one form of $F$. Then $\alpha$ satisfies the Maurer-Cartan equation $d \alpha + \frac{1}{2}[\alpha \wedge \alpha] = 0$~\cite{burstall-pedit}. By decomposing $\alpha$ into its $\goth{h}$ and $\goth{m}$ parts one obtains 
$$
\alpha = \alpha_{\goth{h}} + \alpha_{\goth{m}}, \quad \alpha_{\goth{h}} \in \Gamma ( \goth{h} \otimes T^*\Sigma), \alpha_{\goth{m}} \in \Gamma ( \goth{m} \otimes T^*\Sigma).
$$    
Also according to the decomposition $T^{\C} \Sigma = T' \Sigma \oplus T''\Sigma$,   $\alpha_{\goth{m}}$ decomposes into its $(1,0)$ and $(0,1)$ parts respectively:   
$   \alpha_{\goth{m}} = \alpha'_{\goth{m}}
 + \alpha''_{\goth{m}}$, hence 
\begin{equation}\label{AlphaOneForm1}
\alpha =    \alpha'_{\goth{m}} + \alpha_{\goth{h}} +  \alpha''_{\goth{m}}.
\end{equation}
It is shown (see~\cite{burstall-pedit}) that the harmonic map equation for $\psi $ in terms of $\alpha$ is given by the equation:
\begin{equation}\label{harmonicmapEquation}
 \overline{\partial} \alpha'_{\goth{m}} + [\alpha_{\goth{h}} \wedge \alpha'_{\goth{m}} ] =0.
\end{equation}
There is the following characterization of harmonicity of maps into a symmetric space. Consider the one parameter family of $\goth{g}$-valued one forms  
\begin{equation}\label{AlphaLambdaGenForm}
    \alpha_{\lambda}: = \lambda^{-1} \alpha'_{\goth{m}} + \alpha_{\goth{h}} + \lambda  \alpha''_{\goth{m}}, \quad \lambda \in \bb{S}^1.
\end{equation}
    
\begin{lemma}\label{HarmonicMapsSymmetric}\cite{burstall-pedit}
 $\psi: \Sigma \to G/H$ is harmonic if and only if  
\begin{equation}\label{MaurerCartanlambda}
   d \alpha_{\lambda} + \frac{1}{2}[\alpha_{\lambda} \wedge \alpha_{\lambda}] = 0, \quad \forall \lambda \in \bb{S}^1.
\end{equation}
\end{lemma}
\bigskip

\bigskip

Let $\phi : \Sigma \to Q$ be the Gauss transform of  an immersed  non-isotropic  marginally trapped surface $f : \Sigma \to \bb{S}^4_1$ with non-zero parallel mean curvature vector.  
A frame  $F = (F_0, F_1, F_2, N_1, N_2) \in SO_+(4,1)$  (in column notation) is    {\it adapted}  to  $f$ or $f$-adapted if  $F e_0 =f$ and     
 $F_1, F_2$ span the tangent space of the immersed surface.  Note that if $F$ is  $f$-adapted, then  $N_1, N_2$ are normal sections of $\nu(f)$. \\
The normal bundle $\nu(f)$ is flat  since   $f$ has parallel mean curvature vector~\cite{elghanmi}, thus we can  assume that the normal frame  $\{ N_1, N_2 \} $ in $F$ is positively oriented  and $\nabla^{\bot}$-parallel along $f$.  
Moreover since  $f$  is conformal we can  rotate  within  the tangent plane   $span\{ F_1,  F_2 \}$, if necessary, so that    $f_z = \frac{e^{u}}{\sqrt{2}}(F_1 - i F_2)$.  Let $F: \tilde{\Sigma} \to  SO_+(4,1)$ be an $f$-adapted frame, where $\tilde{\Sigma}$ is the universal covering space of $\Sigma$. Then the structure equations~\eqref{structureEqs} of the immersed surface  $f$   respect to a  coordinate $z$ can be  written as $F_z = F.A$, where    
\begin{equation}\label{matrA}
A =
\begin{pmatrix}
0  &   -\frac{e^u}{\sqrt{2}}  &   i\frac{e^u}{\sqrt{2}}  &    0  & 0\\

\frac{e^u}{\sqrt{2}}& 0  &  i u_z &  - a_1   & a_2\\

-i \frac{e^u}{\sqrt{2}} & -i u_z & 0 & -i  b_1  &  i b_2  \\

0 & a_1  &  i b_1 &       0       &     \sigma        \\

0 & a_2  &  i b_2 &    \sigma     &        0              \\        
\end{pmatrix},
\end{equation}
where the    coefficients in this case are given by    
\begin{equation}\label{ab}
\begin{array}{ll} 
a_1 =  \frac{ e^{u}h-e^{-u} \xi_1}{\sqrt{2}}, &   b_1 =  \frac{ -e^{u} h- e^{-u} \xi_1 }{\sqrt{2}},\\
a_2 =  \frac{ e^{u}h + e^{-u} \xi_2 }{\sqrt{2}}, &   b_2 =   \frac{ - e^u  h+ e^{-u} \xi_2 }{\sqrt{2}}.\\
\end{array}
\end{equation}     
%
%
Defining    $B: = \overline{A}$,   then  $F_{\bar{z}} = FB$ and  the compatibility  among~\eqref{structureEqs}  is just the integrability condition $F_{z \bar{z}} = F_{\bar{z}z}$ which in terms of $A,B$ is given   by  the matrix differential equation $A_{\bar{z}} -  B_z = [A,B]$ encoding   
 Gauss, Codazzi and Ricci's equations~\eqref{gausscodazziricciBasic}.		 In terms of the Maurer-Cartan $\goth{so}(4,1)$-valued one form $\alpha:=F^{-1}dF  =  A dz + B d\bar{z}$ the integrability condition $F_{z \bar{z}} = F_{\bar{z}z}$ is expressed    by the Maurer-Cartan equation $d \alpha + \frac{1}{2} [ \alpha \wedge \alpha] =0$, which is   just the  integrability equation for the existence of an   adapted frame   solving the equation  $F^{-1} dF = \alpha$.  Since $F$ is $f$-adapted then $F$   is  also a  frame for    the Gauss transform $\phi$:
  $$
	\phi =   [ f_z] = [ \frac{e^{u}}{\sqrt{2}}(F_1 - iF_2)]=[F_1 - iF_2] = [F.(e_1-ie_2)] = F.[e_1-ie_2] = F.o 
	$$

We now  decompose    $A = A_{\goth{m}} + A_{\goth{h}} $, and $B = B_{\goth{m}} + B_{\goth{h}} $, where   
\begin{equation}  
A_{\goth{m}} =
\begin{pmatrix}\label{matrAmBm}
0  &   -\frac{e^u}{\sqrt{2}}  &   i\frac{e^u}{\sqrt{2}}  &    0  & 0\\

\frac{e^u}{\sqrt{2}}& 0  &   0 &  - a_1   & a_2\\

-i \frac{e^u}{\sqrt{2}} &  0 & 0 & -i  b_1  &  i b_2  \\

0 & a_1  &  i b_1 &       0       &    0        \\

0 & a_2  &  i b_2 &        0     &        0              \\        
\end{pmatrix}, \quad 
B_{\goth{m}} = \overline{A}_{\goth{m}},
\end{equation}    

\begin{equation}\label{matrAkBk}
A_{\goth{h}} = diag (0, 
(\begin{smallmatrix} 
0 &  i u_{z} \\
-i u_{z} & 0 
\end{smallmatrix}), 
(\begin{smallmatrix} 
0 &  0 \\
0 & 0 
\end{smallmatrix})),  \quad 
B_{\goth{h}} = \overline{A}_{\goth{h}}. 
\end{equation}

Then  
 $A_{\goth{m}}, B_{\goth{m}}$ are   $\goth{m}^{\C}$-valued while       
 $A_{\goth{h}}$ and  $B_{\goth{h}} $ are  $\goth{h}^{\C}$-valued.   Also since  $d \tau_e ([A_{\goth{m}}, B_{\goth{m}} ]) = [A_{\goth{m}}, B_{\goth{m}} ] $,   then      $[A_{\goth{m}}, B_{\goth{m}} ]$ is $\goth{h}^{\C}$-valued. 
Note that  $\alpha'_{\goth{m}}=A_{\goth{m}}dz$,  $\alpha''_{\goth{m}}=B_{\goth{m}}d\bar{z}$ and $\alpha_{\goth{h}}=A_{\goth{h}}dz+ B_{\goth{h}}d\bar{z}$ and so  the harmonic map equation~\eqref{harmonicmapEquation}  for the Gauss transform  $\phi$  becomes  
$$
\Bz A_{\goth{m}} + [B_{\goth{h}}, A_{\goth{m}}] =0.
$$
Here  the family of one forms $\alpha_{\lambda}$~\eqref{AlphaLambdaGenForm}  is given by
\begin{equation}\label{alphalambda}
\alpha_{\lambda} =   \lambda^{-1} A_{\goth{m}}dz + (A_{\goth{h}}dz + B_{\goth{h}}d \bar{z}) + \lambda B_{\goth{m}} d \bar{z}.
\end{equation}

According to Lemma~\ref{HarmonicMapsSymmetric} the Gauss transform $\phi$ is harmonic if and only if    $\alpha_{\lambda}$ satisfies~\eqref{MaurerCartanlambda}. Fixing a point $x_0 \in \tilde{\Sigma}$ and  integrating  for each $\lambda \in \bb{S}^1$     
\begin{equation}\label{dFlambdaODE} 
dF^{\lambda} = F^{\lambda} \alpha_{\lambda}, 
\end{equation}
with initial condition $  F^{\lambda}(x_0) = F(x_0) \in H$, one obtains   a solution  $F^{\lambda} : \tilde{\Sigma} \to  SO_+(4,1)$,   (hence a local solution around any point of  $\Sigma$) which is called an {\it extended frame} normalized at $x_0$. 
  It is possible to choose the constants of integration so that $F^{\lambda}$ depends smoothly on $\lambda \in \bb{S}^1$~\cite{burstall-pedit}. Since $\alpha_{\lambda =1} = \alpha$,  the extended frame satisfies $F^{\lambda=1}(x) = F(x), \forall  x \in \tilde{\Sigma}$.  In column notation,  
$$
F^{\lambda} = (F_0^{\lambda}, F_1^{\lambda}, F_2^{\lambda}, N_1^{\lambda}, N_2^{\lambda}). 
$$
Since  the orthonormal frame  $\{N_1, N_2 \} \subset \Gamma(\nu(f))$ is positively oriented,  an elementary  argument shows that  $\{ N_1^{\lambda}, N_2^{\lambda} \}$ is positively oriented $\forall \lambda \in \bb{S}^1$.\\

Now let    $f^{\lambda} := F_0^{\lambda} =F^{\lambda}e_0$, i.e.  the first column of the extended frame $F^{\lambda}$.  Then  $f^{\lambda}$  is a one parameter deformation  of $f$ since at $\lambda =1$ we recover $f$:  $F^{\lambda=1}e_0 =F.e_0 = f$. We call $f^{\lambda}, \lambda \in \bb{S}^1$ the {\it associated family} of the marginally trapped surface $f$.  
Observe that     
\begin{equation}\label{fzlambda}
f_z^{\lambda} := F_z^{\lambda}e_0 =  F^{\lambda} (\lambda^{-1} A_{\goth{m}} + A_{\goth{h}}) e_0 = \lambda^{-1} \frac{e^{u}}{\sqrt{2}}  F^{\lambda}(e_1 - ie_2),
\end{equation}
hence   $F^{\lambda}$ is  adapted to $f^{\lambda}$. 
From~\eqref{fzlambda}   we compute 
\begin{equation}\label{fzlambdacalculus}
\begin{array}{ll}  
 \la f_z^{\lambda}, f_z^{\lambda} \ra = \la \lambda^{-1} \frac{e^{u}}{\sqrt{2}}  (e_1 - ie_2), \lambda^{-1} \frac{e^{u}}{\sqrt{2}}  (e_1 - ie_2) \ra = 0. & \\
\\
\la f_z^{\lambda}, f_{\bar{z}}^{\lambda} \ra =  \la F^{\lambda} (\lambda^{-1} A_{\goth{m}} + A_{\goth{h}}) e_0, F^{\lambda} (\lambda B_{\goth{m}} + B_{\goth{h}}) e_0 \ra = \\
=\la \lambda^{-1} \frac{e^{u}}{\sqrt{2}}  (e_1 - ie_2), \lambda \frac{e^{u}}{\sqrt{2}}  (e_1 + ie_2) \ra = e^{2u}. 
\end{array}
\end{equation}
Hence      $f^{\lambda}$ is a conformal spacelike immersion  inducing  the same conformal metric  for any $\lambda \in \bb{S}^1$.

Let    $\phi^{\lambda} := F^{\lambda}.o $. Since  $\phi^{ \{\lambda=1\}} = \phi$,    $\phi^{\lambda}$ is a one parameter deformation of $\phi$. The family of maps $\phi^{\lambda}$ is   called   {\it the associated family of} the harmonic Gauss transform  $\phi$~\cite{burstall-pedit}.  
Note that from~\eqref{fzlambda} it follows that   $\phi^{\lambda}$ is the Gauss transform of $f^{\lambda}$:  
$$
\phi^{\lambda}= F^{\lambda}.o = [F^{\lambda}(e_1- i e_2)] = [\lambda f_z^{\lambda}]= [ f_z^{\lambda}]. $$ 

 Hence $\phi^{\lambda}$  takes values in the complex quadric $Q$. Moreover  since   $(\alpha_{\lambda})'_{\goth{m}} = \lambda^{-1} \alpha'_{\goth{m}} $,  $(\alpha_{\lambda})''_{\goth{m}} = \lambda \alpha''_{\goth{m}}$,  and $(\alpha_{\lambda})_{\goth{h}} = \alpha_{\goth{h}}$, it follows that   $\alpha_{\lambda}$ satisfies equation~\eqref{harmonicmapEquation}, thus  
  each  $\phi^{\lambda} : \Sigma \to Q$ is harmonic hence each member $f^{\lambda}$ has parallel mean curvature.  \\

We claim  that $f^{\lambda}$ is marginally trapped for any $\lambda \in \bb{S}^1$. 
Denote by  $\h_{\lambda}$   the mean curvature vector of $f^{\lambda}$. Since $f^{\lambda}$ is conformal and spacelike, it follows that 
\begin{equation}\label{structureOnelambda} 
f^{\lambda}_{z \bar{z}} = -e^{2u} f^{\lambda} + e^{2u} \h_{\lambda}, 
\end{equation}
hence from~\eqref{meancurvature} we obtain   
\begin{equation}\label{meancurvaturelambda}
\h_{\lambda} = e^{-2u}\la f^{\lambda}_{z \bar{z}}, N_1^{\lambda} \ra N_1^{\lambda}- e^{-2u} \la f^{\lambda}_{z \bar{z}}, N_2^{\lambda} \ra N_2^{\lambda}.  
\end{equation}
 
On the other hand  the structure equations of $f^{\lambda}$ are expressed by the matrix equation       
$F_z^{\lambda} = F^{\lambda} (\lambda^{-1} A_{\goth{m}} + A_{\goth{h}})$,  which is equivalent to the system  

\begin{equation}\label{structureeqs2} 
\begin{array}{l}
 f^{\lambda}_z = \frac{1}{\lambda}\frac{e^{u}}{\sqrt{2}} F^{\lambda}_1 - i \frac{1}{\lambda}\frac{e^{u}}{\sqrt{2}} F^{\lambda}_2, \\
 \Pz F^{\lambda}_1= -\frac{1}{\lambda}\frac{e^{u}}{\sqrt{2}} f^{\lambda} -i u_z F^{\lambda}_2 + \frac{1}{\lambda}a_1 N^{\lambda}_1 + \frac{1}{\lambda}a_2 N^{\lambda}_2, \\
\Pz F^{\lambda}_2= i\frac{1}{\lambda}\frac{e^{u}}{\sqrt{2}} f^{\lambda} + i  u_z F^{\lambda}_1 + i\frac{1}{\lambda} b_1 N^{\lambda}_1 + i \frac{1}{\lambda} b_2 N^{\lambda}_2,  \\
\Pz N^{\lambda}_1 = -\frac{1}{\lambda}a_1  F^{\lambda}_1 - i \frac{1}{\lambda}b_1  F^{\lambda}_2+ \sigma N^{\lambda}_1,\\
\Pz N^{\lambda}_2 = \frac{1}{\lambda} a_2  F^{\lambda}_1 + i\frac{1}{\lambda} b_2  F^{\lambda}_2+ \sigma N^{\lambda}_2, \\
\end{array}
\end{equation}
 from which  it follows  
\begin{equation*}
\begin{array}{l}
\la f^{\lambda}_{z \bar{z}}, N_1^{\lambda} \ra =- \la f^{\lambda}_{\bar{z} }, \Pz N_1^{\lambda} \ra =  (a_1 - b_1) \frac{e^{u}}{\sqrt{2}} = e^{2u}h,\\
\la f^{\lambda}_{z \bar{z}}, N_2^{\lambda} \ra = - \la f^{\lambda}_{\bar{z} }, \Pz N_2^{\lambda} \ra =  (b_2 - a_2)  \frac{e^{u}}{\sqrt{2}} = -e^{2u}h.
\end{array}
\end{equation*}
From~\eqref{meancurvaturelambda}  we obtain $\h_{\lambda} = h(N^{\lambda}_1 + N^{\lambda}_2)$, which shows that $f^{\lambda}$ is marginally trapped  for every $\lambda \in \bb{S}^1$, with  $h^{\lambda} = h$. \\

On the other hand since     
 $\xi^{\lambda}_1  = \la f^{\lambda}_{zz}, N^{\lambda}_1 \ra$ and $\xi^{\lambda}_2  =- \la f^{\lambda}_{zz}, N^{\lambda}_2 \ra$, then from~\eqref{structureeqs2} we obtain 
\begin{equation}\label{hopflambda}
\begin{array}{l}
      \xi^{\lambda}_1 = \la f^{\lambda}_{zz}, N_1^{\lambda} \ra= - \la f^{\lambda}_{z}, \Pz N_1^{\lambda} \ra =  \lambda^{-2}\frac{e^{u}}{\sqrt{2}}(a_1 + b_1) =  \lambda^{-2} \xi_1,\\
      \xi^{\lambda}_2 = -\la f^{\lambda}_{zz}, N_2^{\lambda} \ra=  \la f^{\lambda}_{z}, \Pz N_2^{\lambda} \ra = \lambda^{-2}\frac{e^{u}}{\sqrt{2}}(a_2+ b_2) = \lambda^{-2} \xi_2.\\			
\end{array}
\end{equation}
Hence   the $(2,0)$ part of the second fundamental form  of $f_{\lambda}$ is given by,
 \begin{equation}\label{2ndFundformLambda}
\II^{\lambda} (\Pz, \Pz) = \lambda^{-2} \xi_1 N^{\lambda}_1 + \lambda^{-2} \xi_2 N^{\lambda}_2.
\end{equation}
 
From the above expression we see that   $Q_{\lambda} = \lambda^{-2} Q$, where  $Q_{\lambda}=\la f^{\lambda}_{zz}, f^{\lambda}_{zz}\ra dz^4$. Hence  $f^{\lambda}$ is non-isotropic for every $\lambda \in \bb{S}^1$.
 We collect these  facts in the following,

\begin{proposition} 
Let   $f: \Sigma \to \bb{S}^4_1$ be   a non-isotropic conformal marginally trapped immersion with non-zero parallel mean curvature vector and let  $f^{\lambda}$ its  associated family  obtained above which is defined on a simply connected open neighborhhod of each point of $\Sigma$.   Then each member  
 $f^{\lambda}$ is a conformal immersion  inducing the same conformal metric for any $\lambda \in \bb{S}^1$. \\ 
Moreover, $f^{\lambda}$ is a non-isotropic  marginally trapped  surface with   non-zero  parallel mean curvature vector.  

\end{proposition}

\bigskip
Since  $f: \Sigma \to \bb{S}^4_1$ has non-zero parallel mean curvature vector its normal bundle is flat thus by Lemma~\ref{holomorphicHopfdiff} the  Hopf quadratic differential $q = (\xi_1- \xi_2) dz^2$ is holomorphic. Since $f$ is non-isotropic  $q$ is never zero on $\Sigma$, thus  for any point $x \in \Sigma$ there is  a local coordinate $z$ such that $q = c dz^2$, for a non-zero real constant $c$. By~\eqref{ExpressionKappa}     $\kappa$ is real in the same coordinate $z$ and so the null Gauss map $G: \Sigma \to \bb{S}^3$ of $f$ is isothermic.		\\

The conformal invariants $\kappa_{\lambda},s_{\lambda}$ and $\delta$-differential of the 
 associated family $f^{\lambda}$ can be computed as follows. 
	From~\eqref{hopflambda} and~\eqref{2ndFundformLambda} we first obtain  the Hopf differential of $f^{\lambda}$: 
\begin{equation}\label{qLambda}
q_{\lambda} = \lambda^{-2} c dz^2 = \lambda^{-2}q. 
 \end{equation}
Since $h_{\lambda}=h$, from the above expression we obtain 
\begin{equation} \label{deltaLambda}
 \delta_{\lambda} =  h q_{\lambda} = \lambda^{-2} ch  dz^2= \lambda^{-2} \delta.
\end{equation}

In polar  form the Hopf differential $q_{\lambda}$ is given by, 
$q_{\lambda} =   |c|e^{i\theta(\lambda)} dz^2 = \lambda^{-2} |c| e^{i \theta} dz^2$. Thus  $ e^{i\theta(\lambda)} = \lambda^{-2}e^{i\theta}$ so that  if  $\lambda = e^{i \varphi}$ then   
\begin{equation} \label{thetaLambda}
 \theta(\lambda)=\theta-2 \varphi.
\end{equation}
  Since $\lambda$ does not depend on $z$, neither does  $\varphi$ and so    $\theta(\lambda)_z = \theta_z$, and $\theta(\lambda)_{zz} = \theta_{zz}  $. 
Taking this into account and applying~\eqref{schwartzian2} to the conformal invariants $\kappa_{\lambda}, s_{\lambda}$  and the delta  differential  $\delta_{\lambda}$ of $f^{\lambda}$, we obtain,
\begin{equation}\label{schwartzian3}
 c h \lambda^{-2}  =    ((u-i \theta)_{z})^2 + (u-i \theta)_{zz} + \frac{s_{\lambda}}{2}.
\end{equation}
Combining the above equation with~\eqref{schwartzian2} gives  the Schwartzian derivative of $G_{\lambda}$:
\begin{equation}\label{schwartzianLambda}
               s_{\lambda} =  s + 2  (\lambda^{-2} -1) ch. 
\end{equation}

Also from~\eqref{ExpressionKappa2}   $\kappa_{\lambda}$  identifies with  $\frac{e^{u +i \theta(\lambda)}}{\sqrt{2}}$, thus from~\eqref{thetaLambda} we obtain
\begin{equation}\label{kappaLambda}
            \kappa_{\lambda} =  \frac{e^{u +i (\theta- 2 \varphi)}}{\sqrt{2}} = \lambda^{-2} \kappa.
\end{equation}   

A straightforward computation using~\eqref{schwartzian3},~\eqref{schwartzianLambda} and~\eqref{kappaLambda}  shows that $\kappa_{\lambda}, s_{\lambda}, \delta_{\lambda}$ obey the fundamental equation~\eqref{Conf.fund.equation} namely,
$$
         (\kappa_{\lambda})_{\bar{z}\bar{z}} + \frac{\overline{s_{\lambda}}}{2} \kappa_{\lambda} = c h \overline{\lambda^{-2}}  \kappa_{\lambda}, \quad \forall \lambda \in \bb{S}^1,
$$
in which      $ c h \overline{\lambda^{-2}} \kappa_{\lambda} =  c h  \kappa$, hence  it is real valued for every $\lambda \in \bb{S}^1$. In particular    $\kappa_{\lambda}, s_{\lambda}$ obey the conformal Codazzi equation: 
$$
  Im \left( (\kappa_{\lambda})_{\bar{z}\bar{z}} + \frac{\overline{s_{\lambda}}}{2} \kappa_{\lambda} \right) =0, \quad \forall \lambda \in \bb{S}^1.
$$
Since $f$ has parallel mean curvature vector, $\delta$ is holomorphic hence  from~\eqref{schwartzian3},~\eqref{schwartzianLambda},~\eqref{kappaLambda} it easily follows  that    $\kappa_{\lambda}, s_{\lambda}$  obey the conformal Gauss equation: 
$$
\frac{(s_{\lambda})_{\bar{z}}}{2} = 3 (\overline{\kappa_{\lambda}})_z. \kappa_{\lambda}  +  \overline{\kappa_{\lambda}}  (\kappa_{\lambda})_z.
$$	
Since $\lambda$ does not depend on $z$ and  $\delta $ is holomorphic, then
  $\delta_{\lambda} = \lambda^{-2} \delta$  is holomorphic  for any $\lambda \in \bb{S}^1$. 	\\

		We have proved the following

\begin{proposition} Let $f^{\lambda}$  be the associated family of a  non-isotropic marginally trapped   surface  $f: \Sigma \to \bb{S}^4_1$ with  non-zero parallel mean curvature vector. 
Then for any $\lambda \in \bb{S}^1$ the conformal invariants and $\delta$-differential of $f^{\lambda}$ are given by 
\begin{equation}\label{symmetriesG} 
    \kappa_{\lambda} = \lambda^{-2} \kappa, \quad s_{\lambda} =  s + 2  (\lambda^{-2} -1) c h, \quad   \delta_{\lambda} = \lambda^{-2} \delta,
\end{equation}
where $q =c dz^2$ and $\delta = ch dz^2$. \\
Moreover, the system consisting of~\eqref{Conf.fund.equation} and  the conformal Gauss and Codazzi equations ~\eqref{confgausscodazzi} is  invariant under   the  spectral symmetry determined  by~\eqref{symmetriesG}.
\end{proposition}

\noindent

Note that  as  consequence of~\eqref{symmetriesG}   the members of the associated family $f^{\lambda}$  are non-congruent, hence the deformation $f \mapsto f^{\lambda}$ is non-trivial. It also follows  that  the isothermic condition is preserved by the spectral symmetry~\eqref{symmetriesG}: if $\kappa$ is real for some coordinate $z$ then in the new coordinate $w = \frac{1}{\lambda} z $  $\kappa_{\lambda}$ is real since 
  $\kappa_{\lambda}dz^2 = \kappa dw^2$.\\

\begin{remark}
 In~\cite{burstall-pedit-pinkall}  the authors  obtain  the following  slightly different  symmetry for the conformal Gauss and Codazzi equations of a constrained Willmore surface  $\psi : \Sigma \to \bb{S}^3$: 
\begin{equation}\label{symmConstWillmore}
 \kappa_{\lambda} = \lambda \kappa, \quad s_{\lambda} = s + (\lambda^2 -1) \eta, \quad \ \eta_{\lambda} = \lambda^2 \eta,
\end{equation}
where  $\eta dz^2$ is an holomorphic quadratic differential satisfying  $\kappa_{\bar{z}\bar{z}} + \frac{\bar{s}}{2} \kappa= Re(\bar{\eta} \kappa)$.
\end{remark}

\bigskip



\subsection{ The  Calapso-Bianchi associated family of  marginally trapped  surfaces with flat normal bundle}
We construct    an  integrable deformation of non-isotropic marginally trapped surfaces with flat normal bundle which is related to the so-called Calapso-Bianchi T-transform of isothermic surfaces in $\bb{S}^3$~\cite{burstall-pedit-pinkall}. The class of marginally trapped surfaces with flat normal bundle in $\bb{S}^4_1$ includes those with non-zero parallel mean curvature vector and also the spacelike isothermic surfaces introduced  by P. Wang in~\cite{wang}. \\

  Recall that a  conformally immersed  surface $\psi: \Sigma \to \bb{S}^3$ is {\it isothermic} if away from umbilics, it can be
conformally parameterized by its curvature lines.  In terms of its conformal invariants a surface $\psi$ is isothermic if each point in $\Sigma$ has a  coordinate $z$ for which   $\kappa$ is real: $\kappa = \overline{\kappa}$~\cite{burstall-pedit-pinkall},~\cite{xiangma-thesis}. In this case  the conformal Gauss and Codazzi's  equations~\eqref{confgausscodazzi}    away of umbilic points reduce to  
\begin{equation}\label{confGaussCodazzi}
\begin{array}{l}
 s_{\bar{z}} = 4 (\kappa^2)_{z},\\
Im(  \kappa_{\bar{z}\bar{z}} + \frac{1}{2} \bar{s} \kappa) = 0.
\end{array}
\end{equation}
Thus away from umbilic points $\kappa $ is non-zero and so both equations combine into Calapso's equation: $\Delta ( \frac{\kappa_{xy}}{\kappa}) + 8 (\kappa^2)_{xy} =0$. 
The Calapso-Bianchi   T-transform acts on an isothermic surface $\psi : \Sigma \to \bb{S}^3$ by deforming   the schwartzian  $s$  and keeping $\kappa$ unchanged:   
\begin{equation}\label{calapsoT}
  s_t = s+t,  \quad  \kappa_t = \kappa, \quad t \in \R,
\end{equation} 
thus giving rise to the so-called associated family $\psi_t$~\cite{burstall-pedit-pinkall}.\\ 	

Let $f: \Sigma \to \bb{S}^4_1$ be    a  non-isotropic marginally trapped surface  with flat normal bundle. Then by Lemma~\ref{holomorphicHopfdiff}  for every point $x \in \Sigma$ there is a local coordinate $z$ such that $q = c dz^2$ for a non-zero real constant $c$. Thus $\kappa$ is real in the same coordinate and so the null Gauss map $G: \Sigma \to \bb{S}^3$ of $f$ is isothermic. 
Conversely, if  $G$ is isothermic, then $\kappa$ and so $q$ is real in some coordinate $z$. Hence $f$ has flat normal bundle if and only if $q$ is constant i.e. $q=cdz^2$ for some non-zero real constant $c$. \\

 The  structure equations of $f$ read~\eqref{structureEqs} in which  $\xi_1 = \xi +c, \xi_2 =\xi$, $\sigma =0$, where  the positively oriented orthonormal frame  
   $\{ N_1, N_2\}$ is  $\nabla^{\bot}$-parallel.   The compatibility equations~\eqref{gausscodazziricci} reduce in this case to 
    \begin{equation}\label{gausscodazziricci2}
    \begin{array}{ll}
       2 u_{\bar{z}z} = -e^{2u} +e^{-2u} (2c \text{Re}(\xi)+c^2),\\
       \xi_{\bar{z}} =  e^{2u} h_z. \\
			0 = Im((\xi + c)\xi  ), \\ 
    \end{array}
    \end{equation}
		where $2 \xi c + c^2 \neq 0 $ since $f$ is non-isotropic. 
		If $h$ is a non-zero constant, then $f$ has non-zero parallel mean curvature vector field and its null Gauss map $G : \Sigma \to \bb{S}^3$ is isothermic and constrained Willmore.  
	 On the other hand if $h$ is a non-constant function satisfying~\eqref{gausscodazziricci2}, then $f$ has flat normal bundle and non-parallel mean curvature vector field.\\
	
Since our considerations are local we consider an $f$-adapted frame $F= (F_0,F_1, F_2, N_1, N_2)\in SO_+(4,1)$ defined on the universal covering space $\tilde{\Sigma}$.  Then   the structure equations of $f$ read $F_z = FA$, where  the coefficients of the matrix $A$ in~\eqref{matrA}   are given in this case  by  
\begin{equation*}\label{aibi} 
\begin{matrix}
a_1 =  \frac{ e^{-u} (\xi+c)+
e^{u}h }{\sqrt{2}}, &  b_1 = \frac{  e^{-u} (\xi+c) - e^u h}{\sqrt{2}},& \sigma =0,\\
a_2 =  \frac{ e^{-u} \xi+
e^{u}h }{\sqrt{2}}, &  b_2 = \frac{ - e^{-u} \xi - e^u h}{\sqrt{2}}. & \\
\end{matrix}
\end{equation*}  
%


We now  introduce  a one-parameter family of matrices given by  
\begin{equation}
A^t= \begin{pmatrix}
0  &   -\frac{e^{u}}{\sqrt{2}}  &   i\frac{e^{u}}{\sqrt{2}}  &    0  & 0\\

\frac{e^{u}}{\sqrt{2}}& 0  &  i u_z &  - a^t_1   & a^t_2\\

-i \frac{e^{u}}{\sqrt{2}} & -i u_z & 0 & -i  b^t_1  &  i b^t_2  \\

0 & a^t_1  &  i b^t_1 &       0       &     0        \\

0 & a^t_2  &  i b^t_2 &     0     &        0              \\        
\end{pmatrix} , B^t = \overline{A^t} \in \goth{so}(4,1)^{\C}, \quad t \in \R,
\end{equation}
with coefficients   
\begin{equation}\label{aibit} 
\begin{matrix}
a^t_1 =  \frac{ e^{-u} (\xi+c)+
e^{u} h^t }{\sqrt{2}}, &  b^t_1 = \frac{  e^{-u} (\xi+c) - e^{u} h^t}{\sqrt{2}},\\
a^t_2 =  \frac{ e^{-u} \xi+
e^{u} h^t }{\sqrt{2}}, &  b^t_2 = \frac{  -e^{-u} \xi - e^{u} h^t}{\sqrt{2}},\\
\end{matrix}
\end{equation}  
where 

\begin{equation}\label{ansatz}
  h^t := h + \frac{t}{2c},  \quad  c \in \R^{\times},\quad t \in \R.
\end{equation} 
   
Note that for $t=0$ we recover $A$, i.e. $A^{t=0} =A$. 

\begin{lemma}
Define  a   one parameter  family of $\goth{so}(4,1)$-valued one-forms by
\begin{equation}\label{alphat}
\alpha_t := A^t dz + B^t d \bar{z},  \quad t\in \R. 
\end{equation}
Then     $\alpha_t$ coincides with $\alpha$  for  $t=0$ and  it   satisfies the Maurer-Cartan equation 
\begin{equation}\label{maurer-cartan-t}
d \alpha_t + \frac{1}{2} [ \alpha_t \wedge \alpha_t] =0, \quad \forall t \in \R,\end{equation}
if and only if $u, \xi, h$ satisfy~\eqref{gausscodazziricci2}. 
\end{lemma} 
{\bf Proof.}
     Since $B^t = \overline{A^t}$ then $\alpha_t$ is $\goth{so}(4,1)$-valued for every $t \in \R$.
		On the other hand $d \alpha_t + \frac{1}{2} [ \alpha_t \wedge \alpha_t] =0$ is equivalent to $(A^t)_{\bar{z}}-(B^t)_z = [A^t,B^t]$ which in turn is equivalent to 
		\begin{equation*}\label{isothermicGCR}
    \begin{array}{ll}
       2 u_{\bar{z}z} = -e^{2u} +e^{-2u} (2c \text{Re}(\xi)+c^2),\\
       \xi_{\bar{z}} =  e^{2u} (h^t)_z. \\
			0 = Im((\xi + c)\xi  ). \\ 
    \end{array}
    \end{equation*}
Since  $(h^t)_z = h_z$ for any $t \in \R$,  the above system is invariant under the symmetry~\eqref{ansatz}  and  it is equivalent to~\eqref{gausscodazziricci2}.\hfill $\square$\\


Since we work  locally, we may transfer the situation to   the universal covering space $\tilde{\Sigma}$ of $\Sigma$ (note that  the case $\tilde{\Sigma}= \bb{S}^2$ is excluded, otherwise being $q$ holomorphic it would vanish).  Thus we can integrate the Maurer-Cartan equation~\eqref{maurer-cartan-t} on $\tilde{\Sigma}$ for each $t$, obtaining a solution $F^t: \tilde{\Sigma} \to SO_+(4,1)$, which is  unique up to left translation by a constant element in $SO_+(4,1)$.  Thus $F^t$  satisfies 
\begin{equation}\label{structureEqft}
 (F^t)^{-1}d F^t = \alpha_t,  \quad F^{0} =F,  
\end{equation}
since $\alpha_0 = \alpha$. According to~\cite{burstall-pedit},~\cite{ferus-pedit} it is possible to  choose the constants of integration so that  $t \mapsto F^t(x)$ is $C^{\infty}$  for every $x \in \tilde{\Sigma}$. 
Denote   by   $F^t := (F_0^t, F_1^t, F_2^t, N_1^t, N_2^t)$ in column notation. Since   $N^0_2= N_2$ is future pointing, then by continuity $N^t_2$ is future pointing for every $t$. Moreover, since  $\{ N^0_1, N^0_2\}= \{ N_1, N_2\}$ is positively oriented, then an elementary  continuity  argument   shows that    $\{ N_1^t, N_2^t \}$ is positively oriented for every $t \in \R$. \\

Define   $f^t := F^t. e_0$,     the first column of $F^t$, then   
\begin{equation}\label{derivzft}
f^t_z = F^t_z e_0 = F^t A^t. (e_1 - i e_2) = \frac{e^{u}}{\sqrt{2}} F^t (e_1-ie_2),
\end{equation}
from which we compute  
\begin{equation*}
\begin{array}{l}
\la f^t_z, f^t_{\bar{z}} \ra = \frac{e^{2u}}{2} \la F^t  (e_1 - i e_2), F^t  (e_1 + i e_2) \ra = e^{2u}, \\
\la f^t_z, f^t_z \ra = \frac{e^{2u}}{2} \la F^t  (e_1 - i e_2), F^t  (e_1 - i e_2) \ra =0,
\end{array}
\end{equation*}
hence      $f^t$ is a conformal spacelike immersion which     induces the same (conformal) metric  for any $t$. Since  $f^{t=0} =f$, $f^t$ is  a one parameter deformation of $f$.
Also from~\eqref{structureEqft} and~\eqref{derivzft} we obtain 
$$
f^t_{z\bar{z}}=u_{\bar{z}}f^t_z+ \frac{e^{u}}{\sqrt{2}} F^t B(e_1-ie_2), \quad 
f^t_{zz} = u_z f_z + \frac{e^{u}}{\sqrt{2}} F^t A(e_1-ie_2),
$$
 which, from the structure of the matrices $A^t, B^t$, become,
\begin{equation}\label{structureft2}
 f^t_{z\bar{z}} =  -e^{2u}f^t + e^{2u} h^t (N^t_1+N^t_2), \quad 
   f^t_{zz} =  2 u_z f^t_z+ (\xi +c) N^t_1 + \xi N^t_2.
\end{equation}
Hence the mean curvature vector of $f^t$ is given by $\h_t = h^t (N_1^t+ N_2^t)$  and so $f^t$ is marginally trapped. Also from~\eqref{structureft2} we see that 
$$
 \la f^t_{zz}, f^t_{zz} \ra = (\xi +c)^2 - \xi^2= 2 \xi c +c^2 = \la f_{zz}, f_{zz} \ra. \quad \forall t \in R, 
$$
hence $f^t$ is non-isotropic.
On the other hand $F^t$ is adapted to $f^t$ since   $F_z^t = F^t A^t$.  From this equation we extract 
$$
\Pz N^t_1 = -a_1^t F^t_1 - ib^t_1 F^t_2, \quad \Pz N^t_2 = a_2^t F^t_1 + ib^t_2 F^t_2, 
$$
which  shows that $f^t$ has flat normal bundle for every $t$ and that   $\{ N_1^t, N_2^t \}$ is a parallel orthonormal frame  with respect to the normal connection $\nabla^{\bot}_t$ of $\nu(f^t)$. \\

%


Equation~\eqref{Conf.fund.equation} relating the conformal invariants and the $\delta$ differential of  $f$ reads  
\begin{equation}\label{confInvG2}
              \kappa_{\bar{z}\bar{z}} + \frac{\bar{s}}{2} \kappa = ch \kappa, \quad  
   c  \in \R^{\times} , \quad \delta = c h dz^2.
\end{equation}
The deformation family $f^t$ obtained above   is locally defined on $\Sigma$  and  is related to~\eqref{calapsoT} hence  we call $f \mapsto f^t$   {\it the Calapso-Bianchi transformation}  of the marginally trapped surface $f$. \\

Since   $(s_t)_{\bar{z}} = s_{\bar{z}}$,  then by Theorem~\ref{fund.theorem.conf.surf.}  $\kappa, s_t$  determine  a unique (up to Moebius transformations of the sphere) conformal immersed isothermic surface       $G^t: \Sigma \to \bb{S}^3$. Since for $t=0$  we recover $s$ in~\eqref{calapsoT}, $G^t$ is the associated family of $G$ or the T-transform of the isothermic surface $G: \Sigma \to \bb{S}^3$. We claim that $G^t$ is the null Gauss map of $f^t$. In fact,  from~\eqref{structureft2} it follows that $q = c dz^2$ is the Hopf differential of  $f^t$. Since    $f^t$ induce the same conformal metric for all $t$, then  $\theta $ in formula~\eqref{ExpressionKappa2} must be an integer multiple of $2 \pi$, and so   $\kappa = \frac{e^{u}}{\sqrt{2}}$ is the (common) normal Hopf differential of the null Gauss map of all $f^t$. Inserting~\eqref{calapsoT}  into~\eqref{confInvG2} yields, 
\begin{equation}
\kappa_{\bar{z}\bar{z}} + \frac{\overline{s_t}}{2} \kappa = c(h+ \frac{t}{2c}) = ch^t \kappa, \quad  \delta_t = ch^t dz^2, 
\end{equation}
where     $\delta_t = ch^t dz^2$ is just the delta differential of $f^t$. Thus  the above equation is  the  evolution of~\eqref{confInvG2} and so    $\kappa, s_t$ are the conformal invariants  of the null Gauss map of $f^t$.  Thus $G^t$ has  conformal invariants $\kappa, s_t$  and so  it coincides up to a Moebius transformation of $\bb{S}^3$ with the null Gauss map  of $f^t$ which is isothermic since  $\kappa$ is real. \\
The transformation $f \mapsto f^t$ also preserves   marginally trapped surfaces which are isothermic or have parallel second fundamental form. 
For instance if  $f$ is isothermic then  for each $x \in \Sigma$ there is a local coordinate $z$ for which  $\Omega= \xi_1 N_1 dz^2+ \xi_2 N_2 dz^2$ is real valued, that is $\xi_1, \xi_2$ are real valued. Thus by Ricci's equation  $f$ has flat normal bundle and so the Hopf differential $q$ is holomorphic by Lemma~\ref{holomorphicHopfdiff} and so   $q =c dz^2$  for a non-zero real constant $c$, with  $\xi_1 - \xi_2=c$.   Hence the function $\xi$  in~\eqref{gausscodazziricci2} satisfying $\xi_1 = \xi +c, \xi_2 =\xi$  must be also real valued. Thus   from~\eqref{aibit} it follows that the normal vector Hopf differential $\Omega^t$ of $f^t$ is also real valued in the same coordinate $z$, which shows that   $f^t$ is isothermic for any $t \in \R$.  \\
On the other hand if $f$  has non-zero parallel mean curvature vector then it has flat normal bundle by~\cite{elghanmi}. Thus  there is a local positive $\nabla^{\bot}$-parallel orthonormal frame $\{N_1, N_2\} \subset \Gamma(\nu(f))$ such that $0=\nabla^{\bot}_{\Pz} \h = h_z (N_1 + N_2) $, thus  $h$ is constant. Since $h_t$ is defined by~\eqref{ansatz} it satisfies $(h^t)_z = h_z$, then $(h^t)_z = 0$ for all $t \in \R$  which shows  that $f^t$ has (non-zero) parallel mean curvature vector for any $t \in \R$.  
We summarize our  discussion in the following

\begin{theorem}\label{CalapsoBianchitheorem}
 Let $f: \Sigma \to \bb{S}^4_1$ be a non-isotropic conformal marginally trapped immersion  with flat normal bundle.
Let $f^t: \tilde{\Sigma} \to \bb{S}^4_1$ be the Calapso-Bianchi deformation family of $f$ obtained  by integration of~\eqref{structureEqft}. Then on $\Sigma$ each  $f^t$ is  locally defined  conformal non-isotropic marginally trapped immersion with  flat normal bundle whose null Gauss map $G^t$ is isothermic for any $t\in \R$.\\ 
 Moreover, the  transformation $f \mapsto f^t$  preserves   isothermic surfaces and surfaces with  non-zero parallel mean curvature vector.
\end{theorem}

\subsection{ An extended deformation  }

 Non-isotropic marginally trapped conformal immersed surfaces in $\bb{S}^4_1$  with   non-zero parallel mean curvature vector  have flat normal bundle~\cite{elghanmi} and   have   isothermic and constrained  Willmore null Gauss maps into $\bb{S}^3$   by Theorem~\ref{GaussConstWillmore}.  In   the previous section we considered  two different one parameter deformations for such surfaces, namely $f^{\lambda}, \lambda \in \bb{S}^1$ and $f^t, t \in \R$. 
Motivated by~\cite{burstall-pedit-pinkall} we show that it is possible to unify    both deformations by defining an (extended) action of $\C-\{ 0\}$ on the set of non-isotropic marginally trapped surfaces with non-zero parallel mean curvature vector.\\
   
Let $f : \Sigma \to \bb{S}^4_1$ be a non-isotropic conformally immersed marginally trapped surface with non-zero parallel mean curvature vector and   $\kappa, s$ be  the conformal invariants of $f$,  $\delta$-differential $\delta  = c h dz^2$, with $h=const \neq 0$ and  quadratic Hopf differential  $q = cdz^2$, for real constant $c \neq 0$.\\
We   extend the  symmetry~\eqref{symmetriesG}  for  $ \lambda \in \C- \{0\}$ by defining

\begin{equation}\label{symmetriesG2} 
    \kappa_{\lambda} =  |\lambda|^2 \lambda^{-2}  \kappa, \quad s_{\lambda} =  s + 2  (\lambda^{-2} -1) ch , \quad   \delta_{\lambda} = \lambda^{-2} \delta.
\end{equation}
 
Thus for  $|\lambda|=1$ above  we recover~\eqref{symmetriesG}. Moreover, since $ch\kappa $ is real, a straightforward calculation shows  that   $\kappa_{\lambda}, s_{\lambda}, \delta_{\lambda}$ above satisfy ~\eqref{Conf.fund.equation} and the conformal Gauss and Codazzi's equation~\eqref{confgausscodazzi}  for every $\lambda \in \C-\{0\}$. Thus $\kappa_{\lambda}, s_{\lambda}, \delta_{\lambda}$ determine the extended associated familiy  $f^{\lambda}$ which for $|\lambda|=1$ restricts to the associated family obtained in the previous section.   \\

We describe  the deformation~\eqref{symmetriesG2}  of  a non-isotropic marginally trapped torus in $\bb{S}^4_1$ with non-zero parallel mean curvature.  The image of the null Gauss map in this case is an  isothermic constrained Willmore torus in $\bb{S}^3$ and by a result of Richter~\cite{richter} (see also~\cite{burstall-pedit-pinkall})   it can be immersed  as a   surface of constant mean curvature  in some riemannian space form.  \\
Let  $f: \Sigma \to \bb{S}^4_1$ be a conformal non-isotropic marginally trapped immersion with non-zero parallel mean curvature then $Q=\la f_{zz}, f_{zz} \ra dz^4$ is holomorphic and non-zero.  Away from the isolated zeros of $Q$ it is possible to choose a local coordinate $z$ such that $Q = dz^4$, or $\la f_{zz}, f_{zz} \ra=1$. When $\Sigma= T^2$ is a $2$-torus  then $Q$ has no zeros at all (otherwise $Q$ would be identically zero). Thus $Q=dz^4$, where  $z$ is a global coordinate on the universal covering $\C$ of $T^2$ which determines a bi-holomorphism $T^2 \cong \C/\Gamma$ for some lattice $\Gamma_0 \subset \C$. 
  We  choose   a positively oriented orthonormal lorentzian frame $\{ N_1, N_2 \} \subset \Gamma(\nu(f)) $   such that 
$$
f_{zz} = 2 u_z f_z + \cosh(C)N_1 + \sinh(C) N_2,
$$ 
where $C= \rho + i \Theta $ is a  complex function.   The new  positively oriented lorentzian frame $\{ N'_1, N'_2 \}$ given by 
\begin{equation*}
\begin{array}{l}
 N'_1 = \cosh (\rho) N_1 + \sinh (\rho) N_2,\\
N'_2 = \sinh (\rho) N_1 + \cosh (\rho) N_2,\\
\end{array}
\end{equation*}
has  structure function $\sigma'=0$ and so  $\{ N'_1, N'_2 \}$  is  $\nabla^{\bot}$-parallel along $f$. Also since     
$$
f_{zz} = 2 u_z f_z + \cos(\Theta)N'_1 + i \sin(\Theta) N'_2,
$$ 
then  Ricci's equation now becomes $0 = \cos(\Theta) \sin(\Theta)$, of which   $\Theta =0$ is a  solution.  For simplicity we drop the primes and keep denoting by $\{ N_1, N_2 \}$ this new $\nabla^{\bot}$-parallel normal frame. The structure equations of $f$ become    
\begin{equation}\label{structureTorus}
    \begin{array}{l}
    f_{zz} =  2 u_z f_z+   N_1, \\
    f_{\bar{z}z} =  -e^{2u}f + e^{2u} h (N_1+N_2),\\
    \Pz N_1 = -h f_z - e^{-2u}  f_{\bar{z}},\\
    \Pz N_2 = h f_z, 0\neq h=const.
    \end{array}
    \end{equation} 
		with compatibility given by the Sinh-Gordon equation  $2 u_{\bar{z}z} = -e^{2u} +e^{-2u}$, of which        $u: \bb{C} \to \R$ is a  doubly periodic solution with respect to the  lattice $\Gamma_0 \subset \C$.
		Solutions to the Sinh-Gordon equation are obtained  by applying the
finite-gap integration method  from theta
functions  defined on auxiliary hyperelliptic Riemann surfaces which arise from inverse scattering theory~\cite{bobenko}.\\
 Since the mean curvature vector is lightlike and non-zero the codimension of the surface $f$ cannot be reduced. Moreover $ (N_1 + N_2)_z = - e^{-2u}  f_{\bar{z}} $ implies that   $f$ cannot   lie in any singular hypersurface of $\bb{S}^4_1$. 
From~\eqref{structureTorus} the Hopf differential of $f$  is given by $q = dz^2$, hence $\theta$ must be an integer multiple of $2 \pi$ in~\eqref{ExpressionKappa2} and so  $\kappa = \frac{e^{u}}{\sqrt{2}}$ which says that   $G:T^2 \to \bb{S}^3$  is an isothermic and constrained Willmore surface since $h$ is a non-zero constant. The fundamental equation~\eqref{Conf.fund.equation} becomes 
$   \kappa_{\bar{z}\bar{z}} + \frac{\bar{s}}{2} \kappa = h \kappa, \quad \delta = h dz^2$.\\

We see from~\eqref{symmetriesG2}   that  $\kappa_{\lambda} = |\lambda|^2 \lambda^{-2} \frac{e^{u}} {\sqrt{2}}$. Also from~\eqref{kappaVsu} it follows that  the $f^{\lambda}$-induced metric has conformal parameter $u$ for every $\lambda \in \C-\{ 0 \}$, thus  all the surfaces in the extended family have the same induced metric. Using formula~\eqref{ExpressionKappa2} we obtain     
the Hopf quadratic differential of $f^{\lambda}$: 
\begin{equation}\label{hopfLambda}
q_{\lambda} = \lambda^{-2} |\lambda|^2   dz^2.
\end{equation}

Thus since   $
\delta_{\lambda} =   h_{\lambda} q_{\lambda} = \lambda^{-2} \delta= \lambda^{-2} hdz^2$, then the $\delta$-differential of $f^{\lambda}$ is given by    
$\delta_{\lambda} = \lambda^{-2} |\lambda|^2 ( \frac{h}{|\lambda|^2}) dz^2$.   Thus the marginally trapped torus $f^{\lambda}$ has  mean curvature function   $h_{\lambda}= \frac{h}{|\lambda|^2}$  which is a non-zero constant since $\lambda$ does not depend on $z$. Hence $f^{\lambda}$ has non-zero parallel mean curvature vector and  so its null Gauss map $G^{\lambda}$ is  constrained Willmore.   In the  new (rotated) coordinate $w: = \frac{|\lambda|}{\lambda} z$,  \,\,     $\kappa_{\lambda}$ is real respect to $w$   since  $\kappa_{\lambda} dz^2 = \kappa dw^2 $ and $\delta_{\lambda} =  \frac{h}{|\lambda|^2} dw^2$, hence $G^{\lambda}$ is isothermic for every $\lambda \in \C-\{0\}$. \\

Note that  for    $t= 2h (\frac{1}{|\lambda|^2} -1)$  we  recover  the Calapso-Bianchi transformation $f^t$  of the marginally trapped torus $f$.  \\

The structure equations of $f^{\lambda}$ in the extended frame $F^{\lambda} = (f^{\lambda}, f^{\lambda}_z,f^{\lambda}_{\bar{z}}, N^{\lambda}_1, N^{\lambda}_2)$, $\lambda \in \C-\{ 0\}$,  are thus given by     
\begin{equation}\label{structureTorus2}
    \begin{array}{l}
    f^{\lambda}_{zz} =  2 u_z f^{\lambda}_z+ \lambda^{-2} |\lambda|^2  N^{\lambda}_1, \\
    f^{\lambda}_{\bar{z}z} =  -|\lambda|^4 e^{2u}f^{\lambda} + |\lambda|^4 e^{2u} \frac{h}{|\lambda|^2} (N^{\lambda}_1+N^{\lambda}_2),\\
    \Pz N^{\lambda}_1 = -\frac{h}{|\lambda|^2} f^{\lambda}_z - |\lambda|^{-4}e^{-2u}  f^{\lambda}_{\bar{z}},\\
    \Pz N^{\lambda}_2 = \frac{h}{|\lambda|^2} f_z, 
    \end{array}
    \end{equation}

\end{document}